\documentclass[a4paper]{article}
\usepackage[latin1]{inputenc}
\usepackage{amsmath,amssymb}
\usepackage{theorem}
\usepackage{graphics}
\usepackage{epic}
\newtheorem{rem}{Remark}
\newtheorem{bio-rem}[rem]{Biological remark}
\newtheorem{thm}{Theorem}
\newtheorem{prop}{Proposition}
\newtheorem{lemma}{Lemma}
\newtheorem{cor}{Corollary}
\newtheorem{defi}{Definition}
\theoremstyle{break}
\newtheorem{defbr}[defi]{Definition}
\newtheorem{thmbr}[thm]{Theorem}
\newtheorem{lembr}[lemma]{Lemma}

\newtheorem{propbr}[prop]{Proposition}

\title{A Microscopic Interpretation for Adaptive Dynamics Trait
  Substitution Sequence Models}

\author{Nicolas Champagnat\footnote{Weierstrass Institute for Applied
    Analysis and Stochastics \newline Mohrenstrasse 39 \newline 10117
    Berlin, Germany \newline E-Mail: \texttt{champagn@wias-berlin.de}}}

\date{}

\begin{document}

\maketitle

\begin{abstract}
  We consider an interacting particle Markov process for Darwinian
  evolution in an asexual population with non-constant population
  size, involving a linear birth rate, a density-dependent logistic
  death rate, and a probability $\mu$ of mutation at each birth event.
  We introduce a renormalization parameter $K$ scaling the size of the
  population, which leads, when $K\rightarrow+\infty$, to a
  deterministic dynamics for the density of individuals holding a
  given trait. By combining in a non-standard way the limits of large
  population ($K\rightarrow+\infty$) and of small mutations
  ($\mu\rightarrow 0$), we prove that a time scales separation between
  the birth and death events and the mutation events occurs and that
  the interacting particle microscopic process converges for finite
  dimensional distributions to the biological model of evolution known
  as the ``monomorphic trait substitution sequence'' model of adaptive
  dynamics, which describes the Darwinian evolution in an asexual
  population as a Markov jump process in the trait space.
\end{abstract}

\medskip \emph{Keywords:} measure-valued process; interacting particle
system; mutation-selection processes; Darwinian evolution; trait
substitution sequence; adaptive dynamics; finite dimensional
distributions convergence; time scale separation; stochastic
domination; branching processes; large deviations.

\medskip
\emph{AMS subject classification:} 60F99; 60K35; 92D15.

\section{Introduction and main results}
\label{sec:intro}

We will study in this article the link between two biological models
of Darwinian evolution in an asexual population. The first one is a
system of interacting particles modeling evolution at the
\emph{individual} level, referred below as the \emph{microscopic
  model}, which has been already proposed and studied in Bolker and
Pacala (1997, 1999), Dieckmann and Law (2000), Law et al.\ (2003) and
Fournier and M\'el\'eard (2004) either as a model of Darwinian
evolution or as a model of dispersal in a spatially structured
population. This model involves a finite population with non-constant
population size, in which each individual's birth and death events are
described. Each individual's ability to survive and reproduce is
characterized by a finite number of phenotypic traits (e.g.\ body
size, rate of food intake, age at maturity), or simply \emph{traits}.
The birth rate of an individual depends on its phenotype, and its
death rate depends on the distribution of phenotypes in the population
and involves a competition kernel of logistic type. A mutation may
occur at each birth event.

The second model describes the evolution at the \emph{population}
level as a jump Markov process in the space of phenotypic traits
characterizing individuals. It is called ``trait substitution
sequence'' (Metz et al., 1996), and referred below as the \emph{TSS
  model}. In this model, the population is \emph{monomorphic} at each
time (i.e.\ composed of individuals holding the \emph{same} trait
value), and the evolution proceeds by a sequence of appearance of new
mutant traits, which invade the population and replace, after a short
competition, the previous dominant trait. The TSS model belongs to the
recent biological theory of evolution called \emph{adaptive dynamics}
(Hofbauer and Sigmund, 1990; Marrow et al., 1992; Metz et al., 1992),
and has been introduced by Metz et al.\ (1996) and Dieckmann and Law
(1996) and mathematically studied in Champagnat et al.\ (2001). The
theory of adaptive dynamics investigates the effects of the ecological
aspects of population dynamics on the evolutionary process, and thus
describes the population on the phenotypic level, instead of the
genotypic level.  The TSS model is one of the fundamental models of
this theory. It has revealed a powerful tool for understanding various
evolutionary phenomena, such as polymorphism (stable coexistence of
different traits, cf.\ Metz et al., 1996) or evolutionary branching
(evolution of a monomorphic population to a polymorphic one that may
lead to speciation, Dieckmann and Doebeli, 1999) and is the basis of
other biological models, such as the ``canonical equation of adaptive
dynamics'' (Dieckmann and Law, 1996; Champagnat et al., 2001).

The heuristics leading to the TSS model (cf.\ Metz et al., 1996 and
Dieckmann and Law, 1996) are based on the biological assumptions of
large population and rare mutations, and on another assumption stating
that no two different types of individuals can coexist on a long time
scale: the competition eliminates one of them.  In spite of this
heuristic, this model still lacks a firm mathematical basis.

We propose to prove in this article a convergence result of the
microscopic model to the TSS model when the parameters are normalized
in a non-standard way, leading to a \emph{time scales separation}. Our
limit combines a \emph{large population} asymptotic with a \emph{rare
  mutations} asymptotic. It will appear that this convergence holds
only for finite dimensional distributions, and not for the Skorohod
topology, for reasons that are linked to the time scale separation.
For these reasons, and because we have to combine two limits
simultaneously (large population and rare mutations), this result is
different from classical time scale separation results (averaging
principle, cf.\ Freidlin and Wentzell, 1984). The proof requires
original methods, based on comparison, convergence and large deviation
results on branching processes and logistic Markov birth and death
processes. Our convergence result provides a mathematical
justification of the TSS model and of the biological heuristic on
which it is based, and gives precise conditions on the scalings of the
biological parameters in the microscopic model required for the time
scales separation to hold.

In Section~\ref{sec:models}, we describe precisely the microscopic
model and the TSS model, and we state our main results.  Our proof is
based on a careful study of the behavior of the population before the
first mutation, and of the competition phase between the mutant trait
and the original trait, taking place just after the first mutation.
We will give an outline of the proof and of the methods in
Section~\ref{sec:outline}, as well as some notations used throughout
the paper. Section~\ref{sec:birth-death-proc} gives comparison results
and large deviation results on birth and death processes
(Sections~\ref{sec:compar} and~\ref{sec:exit-dom}), and several
results on branching processes (Section~\ref{sec:br-proc}). Based on
these properties, the proof of the convergence of the microscopic
model to the TSS model is given in Section~\ref{sec:proof}.

\section{Models and main results}
\label{sec:models}

Let us first describe the microscopic model. In a population,
Darwinian evolution acts on a set of phenotypes, or \emph{traits},
characterizing each individual's ability to survive and reproduce. We
consider a finite number of quantitative traits in an asexual
population (clonal reproduction), and we assume that the trait space
${\cal X}$ is a compact subset of $\mathbb{R}^l$ ($l\geq 1$).

The microscopic model involves the three basic mechanisms of Darwinian
evolution: \emph{heredity}, which transmits traits to new offsprings,
\emph{mutation}, driving a variation in the trait values in the
population, and \emph{selection} between these different trait values.
The selection process, and thus a proper definition of the selective
ability of a trait, or \emph{fitness} (cf.\ Metz et al., 1992), should
(and will) be the consequence of interactions between individuals in
the population and of the competition for limited resources or area,
modeled as follows.

For any $x,y\in{\cal X}$, we introduce the following biological
parameters
\begin{description}
\item[$b(x)\in\mathbb{R}_+$] is the rate of birth from an individual
  holding trait $x$.
\item[$d(x)\in\mathbb{R}_+$] is the rate of ``natural'' death for an
  individual holding trait $x$.
\item[$\alpha(x,y)\in\mathbb{R}_+$] is competition kernel
  representing the pressure felt by an individual holding trait $x$
  from an individual holding trait $y$.
\item[\mbox{$\mu(x)\in[0,1]$}] is the probability that a mutation
  occurs in a birth from an individual with trait $x$.
\item[$m(x,dh)$] is the law of $h=y-x$, where the mutant trait $y$ is
  born from an individual with trait $x$. It is a probability measure
  on $\mathbb{R}^l$, and since $y$ must belong to the trait space ${\cal
    X}$, the support of $m(x,\cdot)$ is a subset of
  \begin{equation*}
    {\cal X}-x=\{ y-x:y\in{\cal X}\}.
  \end{equation*}
\item[$K\in\mathbb{N}$] is a parameter rescaling the competition
  kernel $\alpha(\cdot,\cdot)$. Biologically, $K$ can be interpreted
  as scaling the resources or area available, and is related to the
  biological concept of ``carrying capacity''. It is also called
  ``system size'' by Metz et al. (1996). As will appear later, this
  parameter is linked to the size of the population: large $K$ means a
  large population (provided that the initial condition is
  proportional to $K$).
\item[\mbox{$u_K\in[0,1]$}] is a parameter depending on $K$ rescaling
  the probability of mutation $\mu(\cdot)$. Small $u_K$ means rare
  mutations.
\end{description}

Let us also introduce the following notations, used throughout this
paper:
\begin{gather}
  \label{eq:def-bar-n}
  \bar{n}_x=\frac{b(x)-d(x)}{\alpha(x,x)}, \\
  \beta(x)=\mu(x)b(x)\bar{n}_x \label{eq:def-beta} \\
  \mbox{and}\quad f(y,x)=b(y)-d(y)-\alpha(y,x)\bar{n}_x. \label{eq:def-fitn}
\end{gather}
As will appear below, $\bar{n}_x$ can be interpreted as the
equilibrium density of a monomorphic population when there is no
mutation, $\beta(x)$ as the mutation rate in this population, and
$f(y,x)$ as the fitness of a mutant individual with trait $y$ in this
population.

We consider, at any time $t\geq 0$, a finite number $N_t$ of
individuals, each of them holding a trait value in ${\cal X}$. Let us
denote by $x_1,\ldots,x_{N_t}$ the trait values of these individuals.
The state of the population at time $t\geq 0$, rescaled by $K$, can be
described by the finite point measure on ${\cal X}$
\begin{equation}
  \label{eq:nu_t}
  \nu^K_t=\frac{1}{K}\sum_{i=1}^{N_t}\delta_{x_i},
\end{equation}
where $\delta_x$ is the Dirac measure at $x$. Let ${\cal M}_F$ denote
the set of finite nonnegative measures on ${\cal X}$, and define
\begin{equation*}
  {\cal M}^K=\left\{\frac{1}{K}\sum_{i=1}^n\delta_{x_i}:n\geq 0,\
  x_1,\ldots,x_n\in{\cal X}\right\},
\end{equation*}

An individual holding trait $x$ in the population $\nu^K_t$ gives
birth to another individual with rate $b(x)$ and dies with rate
\begin{equation*}
  d(x)+\int\alpha(x,y)\nu^K_t(dy)=d(x)
  +\frac{1}{K}\sum_{i=1}^{N_t}\alpha(x,x_i).
\end{equation*}
The parameter $K$ scales the strength of competition, thus allowing
the coexistence of more individuals in the population.

A newborn holds the same trait value as its progenitor's with
probability $1-u_K\mu(x)$, and with probability $u_K\mu(x)$, the
newborn is a mutant whose trait value $y$ is chosen according to
$y=x+h$, where $h$ is a random variable with law $m(x,dh)$.

In other words, the process $(\nu^K_t,t\geq 0)$ is a ${\cal
  M}^K$-valued Markov process with infinitesimal generator defined for
any bounded measurable functions $\phi$ from ${\cal M}^K$ to
$\mathbb{R}$ by
\begin{align}
  L^K\phi(\nu) & =\int_{\cal
    X}\left(\phi\left(\nu+\frac{\delta_x}{K}\right)
    -\phi(\nu)\right)(1-u_K\mu(x))b(x)K\nu(dx) \notag \\
  & +\int_{\cal X}\int_{\mathbb{R}^l}
  \left(\phi\left(\nu+\frac{\delta_{x+h}}{K}\right)-\phi(\nu)\right)
  u_K\mu(x)b(x)m(x,dh)K\nu(dx) \notag \\
  & +\int_{\cal X}\left(\phi\left(\nu-\frac{\delta_x}{K}
    \right)-\phi(\nu)\right)\left(d(x)+\int_{\cal
      X}\alpha(x,y)\nu(dy)\right)K\nu(dx).
  \label{eq:generator-renormalized-IPS}
\end{align}
When the measure $\nu$ has the form~(\ref{eq:nu_t}), the integrals
with respect to $K\nu(dx)$ in~(\ref{eq:generator-renormalized-IPS})
correspond to sums over all individual in the population.  The first
term (linear) describes the births without mutation, the second term
(linear) describes the births with mutation, and the third term
(non-linear) describes the deaths by oldness or competition. This
logistic density-dependence models the competition in the population,
and hence drives the selection process.

Let us denote by~(A) the following three assumptions
\begin{description}
\item[\textmd{(A1)}] $b$, $d$ and $\alpha$ are measurable functions,
  and there exist $\bar{b},\bar{d},\bar{\alpha}<+\infty$ such that
  \begin{equation*}
    b(\cdot)\leq\bar{b},\quad d(\cdot)\leq\bar{d}
    \quad\mbox{and}\quad\alpha(\cdot,\cdot)\leq\bar{\alpha}.
  \end{equation*}
\item[\textmd{(A2)}] $m(x,dh)$ is absolutely continuous with respect
  to the Lebesgue measure on $\mathbb{R}^l$ with density $m(x,h)$, and
  there exists a function $\bar{m}:\mathbb{R}^l\rightarrow\mathbb{R}_+$ such
  that $m(x,h)\leq \bar{m}(h)$ for any $x\in{\cal X}$ and
  $h\in\mathbb{R}^l$, and $\int \bar{m}(h)dh<\infty$.
\item[\textmd{(A3)}] $\mu(x)>0$ and $b(x)-d(x)>0$ for any $x\in{\cal
    X}$, and there exists $\underline{\alpha}>0$ such
  that
  \begin{equation*}
    \underline{\alpha}\leq\alpha(\cdot,\cdot).
  \end{equation*}
\end{description}

For fixed $K$, under~(A1) and~(A2) and assuming that
$\mathbf{E}(\langle\nu^K_0,\mathbf{1}\rangle)<\infty$ (where
$\langle\nu,f\rangle$ denotes the integral of the measurable function
$f$ with respect to the measure $\nu$), the existence and uniqueness
in law of a process with infinitesimal generator $L^K$ has been proved
by Fournier and M\'el\'eard (2003). When $K\rightarrow+\infty$, they
also proved, under more restrictive assumptions and assuming the
convergence ot the initial condition, the convergence on
$\mathbb{D}(\mathbb{R}_+,{\cal M}_F)$ of the process $\nu^K$ to a
deterministic process solution to a non-linear integro-differential
equation. We will only use particular cases of their result, stated in
the next section, that can be proved under assumptions~(A1) and~(A2).

The biological assumption of large population corresponds to the limit
$K\rightarrow+\infty$, and the assumption of rare mutations to
$u_K\rightarrow 0$. As mentionned in the introduction, the biological
heuristics suggest another assumption: the impossibility of
coexistence of two different traits on a long time scale. As will
appear in Proposition~\ref{prop:LV} in the next section, this
assumption can be stated mathematically as follows:
\begin{description}
\item[\textmd{(B)}] Given any $x\in{\cal X}$, Lebesgue almost any
  $y\in{\cal X}$ satisfies one of the two following conditions:
  \begin{align}
    \mbox{either}\quad &
    (b(y)-d(y))\alpha(x,x)-(b(x)-d(x))\alpha(y,x)<0, \label{eq:hyp-B1} \\
    \mbox{or}\quad & \left\{
      \begin{array}{l}
        (b(y)-d(y))\alpha(x,x)-(b(x)-d(x))\alpha(y,x)>0, \\
        (b(x)-d(x))\alpha(y,y)-(b(y)-d(y))\alpha(x,y)<0.
      \end{array} \right. \label{eq:hyp-B2}
  \end{align}
\end{description}
Before coming back to this assumption in the next section, let us only
observe that condition~(\ref{eq:hyp-B1}) is equivalent to $f(y,x)<0$
and condition~(\ref{eq:hyp-B2}) to $f(y,x)>0$ and $f(x,y)<0$.

The TSS model of evolution that we obtain from the microscopic model
is a Markov jump process in the trait space ${\cal X}$ with
infinitesimal generator given, for any bounded measurable function
$\varphi$ from ${\cal X}$ to $\mathbb{R}$, by
\begin{equation}
  \label{eq:generator-TSS}
  A\varphi(x)=\int_{\mathbb{R}^l}(\varphi(x+h)-\varphi(x))
  \beta(x)\frac{[f(x+h,x)]_+}{b(x+h)}m(x,h)dh,
\end{equation}
where $[a]_+$ denotes the positive part of $a\in\mathbb{R}$, and where
$\beta(x)$ and $f(y,x)$ are defined in~(\ref{eq:def-beta})
and~(\ref{eq:def-fitn}).  The existence and uniqueness in law of a
process generated by $A$ holds as soon as $\beta(x)[f(y,x)]_+/b(y)$ is
bounded (see e.g.\ Ethier and Kurtz, 1986), which is true under
assumption~(A) ($[f(y,x)]_+/b(y)\leq 1$). The biological
interpretation of the fonction $f$ as a \emph{fitness} function
becomes natural in view of this generator: because of the positive
part function $[\cdot]_+$ in~(\ref{eq:generator-TSS}), the TSS process
can only jump from a trait $x$ to the traits $x+h$ such that
$f(x+h,x)>0$. Therefore, the function $f(y,x)$ measures the selective
ability of trait $y$ in a population made of individuals with trait
$x$ (see Metz et al., 1992, 1996).

Our main result is:
\begin{thm}
  \label{thm:IPS-TSS}
  Assume~(A) and~(B). Fix a sequence $(u_K)_{K\in\mathbb{N}}$ in
  $[0,1]^{\mathbb{N}}$ such that
  \begin{equation}
    \label{eq:cond-mu-K}
    \forall V>0,\quad\exp(-VK)\ll u_K\ll\frac{1}{K\log K}
  \end{equation}
  (where $f(K)\ll g(K)$ means that $f(K)/g(K)\rightarrow 0$ when
  $K\rightarrow\infty$). Fix also $x\in{\cal X}$, $\gamma>0$ and a
  sequence of $\mathbb{N}$-valued random variables
  $(\gamma_K)_{K\in\mathbb{N}}$, such that
  $(\gamma_K/K)_{K\in\mathbb{N}}$ converges in law to $\gamma$ and is
  bounded in $\mathbb{L}^1$. Consider the process $(\nu^K_t,t\geq 0)$
  generated by~(\ref{eq:generator-renormalized-IPS}) with initial
  state $(\gamma_K/K)\delta_x$. Then, for any $n\geq 1$,
  $\varepsilon>0$ and $0<t_1<t_2<\ldots<t_n<\infty$, and for any
  measurable subsets $\Gamma_1,\ldots,\Gamma_n$ of ${\cal X}$,
  \begin{multline}
    \label{eq:def-cvgce}
    \lim_{K\rightarrow+\infty}\mathbf{P}\bigl(\forall i\in\{1,\ldots,n\},\
    \exists x_i\in\Gamma_i:\mbox{\textnormal{Supp}}(\nu^K_{t_i/K
    u_K})=\{x_i\} \\ \mbox{and\ }|\langle
    \nu^K_{t_i/K u_K},\mathbf{1}\rangle-\bar{n}_{x_i}|<\varepsilon\bigr)
    =\mathbf{P}(\forall i\in\{1,\ldots,n\},\ X_{t_i}\in\Gamma_i)
  \end{multline}
  where for any $\nu\in{\cal M}_F$, $\mbox{\textnormal{Supp}}(\nu)$ is
  the support of $\nu$ and $(X_t,t\geq 0)$ is the TSS process
  generated by~(\ref{eq:generator-TSS}) with initial state $x$.
\end{thm}

\begin{rem}
  \label{rem:initial-cond}
  The time scale $1/Ku_K$ of Theorem~\ref{thm:IPS-TSS} is the time
  scale of the mutation events for the process $\nu^K$ (the population
  size is proportional to $K$ and the individual mutation rate is
  proportional to $u_K$). Assumption~(\ref{eq:cond-mu-K}) is the
  condition leading to the correct time scales separation between the
  mutation events and the birth and death events. The
  limit~(\ref{eq:def-cvgce}) means that, when this time scales
  separation occurs, the population is monomorphic at any time with
  high probability, and that the transition periods corresponding to
  the invasion of a mutant trait in the resident population and the
  ensuing competition are infinitesimal on this mutation time scale.
  Observe also that this convergence result holds only for monomorphic
  initial conditions. We will make some comments on more general
  initial conditions in the next section.
\end{rem}

\begin{cor}
  \label{cor:IPS-TSS}
  Assume additionally in Theorem~\ref{thm:IPS-TSS} that
  $(\gamma_K/K)_{K\in\mathbb{N}}$ is bounded in $\mathbb{L}^p$ for
  some $p>1$. Then the process $(\nu^K_{t/Ku_K},t\geq 0)$ converges
  when $K\rightarrow +\infty$, in the sense of the finite dimensional
  distributions for the topology on ${\cal M}_F$ induced by the
  functions $\nu\mapsto\langle\nu,f\rangle$ with $f$ bounded and
  measurable on ${\cal X}$, to the process $(Y_t,t\geq 0)$ defined by
  \begin{equation*}
    Y_t=\left\{\begin{array}{ll}
        \gamma\delta_x & \mbox{if\ }t=0 \\
        \bar{n}_{X_t}\delta_{X_t} & \mbox{if\ }t>0.
      \end{array}\right.
  \end{equation*}
\end{cor}

This corollary follows from the following long time moment estimates.
\begin{lemma}
  \label{lem:moment}
  Assume~(A) and that $\:\sup_{K\geq 1} \mathbf{E}(\langle
  \nu^K_0,1\rangle^p)<+\infty$ for some $p\geq 1$, then
  \begin{equation*}
    \sup_{K\geq 1}\:\sup_{t\geq 0}\mathbf{E}\big(\langle
    \nu^K_t,\mathbf{1}\rangle^p\big)<+\infty,
  \end{equation*}
  and therefore, if $p>1$, the family of random variables $\{\langle
  \nu^K_t,\mathbf{1}\rangle\}_{\{K\geq 1,\:t\geq 0\}}$ is uniformly
  integrable.
\end{lemma}

\paragraph{Proof of Corollary~\ref{cor:IPS-TSS}}
Let $\Gamma$ be a measurable subset of ${\cal X}$. Let us prove that
\begin{equation}
  \label{eq:pf-cor}
  \lim_{K\rightarrow+\infty}
  \mathbf{E}(\langle\nu^K_{t/Ku_K},\mathbf{1}_{\Gamma}\rangle)
  =\mathbf{E}(\bar{n}_{X_t}\mathbf{1}_{\{X_t\in\Gamma}\}).
\end{equation}
Fix $\varepsilon>0$, and observe that
$\bar{n}_x\in[0,\bar{b}/\underline{\alpha}]$. Write
$[0,\bar{b}/\underline{\alpha}]\subset\cup_{i=1}^{q}I_i$, where $q$ is
the first integer greater than
$\bar{b}/\varepsilon\underline{\alpha}$, and
$I_i=[(i-1)\varepsilon,i\varepsilon[$. Define $\Gamma_i=\{x\in{\cal
  X}:\bar{n}_x\in I_i\}$ for $1\leq i\leq q$, and
apply~(\ref{eq:def-cvgce}) to the sets
$\Gamma\cap\Gamma_1,\ldots,\Gamma\cap\Gamma_q$ with $n=1$, $t_1=t$ and
the constant $\varepsilon$ above. Then, by Lemma~\ref{lem:moment},
there exists a constant $C>0$ such that
\begin{align*}
  \limsup_{K\rightarrow+\infty}
  \mathbf{E}(\langle\nu^K_{t/Ku_K},\mathbf{1}_{\Gamma}\rangle)
  & \leq\limsup_{K\rightarrow+\infty}
  \mathbf{E}(\langle\nu^K_{t/Ku_K},\mathbf{1}_{\Gamma}\rangle
  \mathbf{1}_{\{\langle\nu^K_{t/Ku_K},\mathbf{1}\rangle\leq C\}})+\varepsilon \\
  & \leq\sum_{i=1}^q\limsup_{K\rightarrow+\infty}
  \mathbf{E}(\langle\nu^K_{t/Ku_K},\mathbf{1}_{\Gamma\cap\Gamma_i}\rangle
  \mathbf{1}_{\{\langle\nu^K_{t/Ku_K},\mathbf{1}\rangle\leq C\}})+\varepsilon \\
  & \leq\sum_{i=1}^q(i+1)\varepsilon\mathbf{P}(X_t\in\Gamma\cap\Gamma_i)
  +\varepsilon \\
  & \leq\sum_{i=1}^q\bigl(
  \mathbf{E}(\bar{n}_{X_t}\mathbf{1}_{\{X_t\in\Gamma\cap\Gamma_i\}})
  +2\varepsilon\mathbf{P}(X_t\in\Gamma_i)\bigr)+\varepsilon \\
  & \leq\mathbf{E}(\bar{n}_{X_t}\mathbf{1}_{\{X_t\in\Gamma\}})+3\varepsilon.
\end{align*}
A similar estimate for the \emph{lim\:inf} ends the proof
of~(\ref{eq:pf-cor}), which implies the convergence of one-dimensional
laws for the required topology.

The same method gives easily the required limit when we consider a
finite number of times $t_1,\ldots,t_n$.\hfill$\Box$
\medskip

As suggested by the fact that the limit process $Y$ is not continuous
at $0^+$, it is not possible to obtain the convergence in law for the
Skorohod topology on $\mathbb{D}([0,T],{\cal M}_F)$. More generally,
we can prove:

\begin{prop}
  \label{prop:non-cvgce}
  For any $s<t$, the convergence of $\nu^K_{\cdot/Ku_K}$ to $Y$ in
  Corollary~\ref{cor:IPS-TSS} does not hold for the Skorohod topology
  on $\mathbb{D}([s,t],{\cal M}_F)$, for any topology on ${\cal M}_F$
  such that the total mass function
  $\nu\mapsto\langle\nu,\mathbf{1}\rangle$ is continuous.
\end{prop}

\paragraph{Proof of Proposition~\ref{prop:non-cvgce}}
Assume the converse. Then, for some $s<t$, the total mass
$N^K_t=\langle\nu^K_{t/Ku_K},\mathbf{1}\rangle$ converges for the
Skorohod topology on $\mathbb{D}([s,t],\mathbb{R}_+)$ to the total
mass of the process $Y$. In particular, by Ascoli's theorem for
c\`adl\`ag processes (cf.\ Billingsley, 1968), for any $\varepsilon>0$
and $\eta>0$, there exists $\delta>0$ such that
\begin{equation*}
  \limsup_{K\rightarrow+\infty}\mathbf{P}(\omega'(N^K,\delta)>\eta)\leq\varepsilon,
\end{equation*}
where the modulus of continuity $\omega'$ is defined by
\begin{equation*}
  \omega'(\varphi,\delta):=\inf\left\{\max_{i=0,\ldots,r-1}
  \omega(\varphi,[t_i,t_{i+1}))\right\}
\end{equation*}
where the infimum is taken over all $r\in\mathbb{N}$ and all the
finite partitions $s=t_0<t_1<\ldots<t_r=t$ of $[s,t]$ such that
$t_{i+1}-t_i>\delta$ for any $i\in\{0,\ldots,r-1\}$, and where
$\omega(\varphi,I):=\sup_{x,y\in I}|\varphi(x)-\varphi(y)|$ for any
interval $I$.

Now, for any function $\varphi\in\mathbb{D}([s,t],\mathbb{R})$,
$\omega(\varphi,\delta)\leq 2\omega'(\varphi,\delta)+\sup_{x\in[s,t]}
|\varphi(x)-\varphi(x-)|$ (cf.\ Billingsley, 1968), where
$\omega(\varphi,\delta):=\sup_{x,y\in[s,t],\ 
  |x-y|\leq\delta}|\varphi(x)-\varphi(y)|$, and for any $K\geq 1$,
$\sup_{x\in[s,t]}|N^K_x-N^K_{x-}|=1/K$. Therefore, for any
$\varepsilon>0$ and $\eta>0$, there exists $\delta>0$ such that
\begin{equation*}
  \limsup_{K\geq 1}\mathbf{P}(\omega(N^K,\delta)>\eta)\leq\varepsilon.
\end{equation*}
This implies that the sequence $(N^K)_K$ is actually C-tight (cf.\ 
Billingsley, 1968) and that its limit is necessarily continuous, which
is not true for $\langle Y_t,\mathbf{1}\rangle$.\hfill$\Box$



\section{Notations and outline of the proof of Theorem~\ref{thm:IPS-TSS}}
\label{sec:outline}

We start with some definitions needed to explain the idea of the proof
of Theorem~\ref{thm:IPS-TSS} and the precise meaning of
assumption~(B).

\begin{defbr}
  \label{def:birth-death-proc}
  \begin{description}
  \item[\textmd{(a)}] For any $K\geq 1$, $b,d,\alpha\geq 0$ and for any
    $\mathbb{N}/K$-valued random variable $z$, we will denote by
    $\mathbf{P}^K(b,d,\alpha,z)$ the law of the $\mathbb{N}/K$-valued
    Markov birth and death process with initial state $z$ and with transition
    rates
    \begin{equation*}
      \begin{array}{ll}
        ib & \mbox{from\ }i/K\mbox{\ to\ }(i+1)/K, \\
        i(d+\alpha i/K) & \mbox{from\ }i/K\mbox{\ to\ }(i-1)/K.
      \end{array}
    \end{equation*}
  \item[\textmd{(b)}] For any $K\geq 1$, $b_k,d_k,\alpha_{kl}\geq 0$ with
    $k,l\in\{1,2\}$, and for any $\mathbb{N}/K$-valued random
    variables $z_1$ and $z_2$, we will denote by
    \begin{equation*}
      \mathbf{Q}^K(b_1,b_2,d_1,d_2,\alpha_{11},\alpha_{12},
      \alpha_{21},\alpha_{22},z_1,z_2)
    \end{equation*}
    the law of the $(\mathbb{N}/K)^2$-valued Markov birth and death with
    initial state $(z_1,z_2)$ and with transition rates
    \begin{equation*}
      \begin{array}{ll}
        ib_1 & \mbox{from\ }(i/K,j/K)\mbox{\ to\ }((i+1)/K,j/K), \\
        jb_2 & \mbox{from $(i/K,j/K)$ to $(i/K,(j+1)/K)$,} \\
        i(d_1+\alpha_{11}i/K+\alpha_{12}j/K) & \mbox{from $(i/K,j/K)$ to
          $((i-1)/K,j/K)$,} \\
        j(d_2+\alpha_{21}i/K+\alpha_{22}j/K) & \mbox{from
          $(i/K,j/K)$ to $(i/K,(j-1)/K)$.}
      \end{array}
    \end{equation*}
  \end{description}
\end{defbr}

These two Markov processes have absorbing states at $0$ and $(0,0)$,
respectively. Observe also that, when $\alpha=0$, the Markov process
of point~(a) is a continuous-time binary branching process divided by
$K$.

Fix $x$ and $y$ in ${\cal X}$. The proof of the following two results
can be found in Chap.~11 of Ethier and Kurtz (1986).
\begin{propbr}
  \label{prop:cv-mono-di}
  \begin{description}
  \item[\textmd{(a)}] Assume $\mu\equiv 0$ and
    $\nu^K_0=N^K_x(0)\delta_x$. Then, for any $t\geq 0$,
    $\nu^K_t=N^K_x(t)\delta_x$, where $N^K_x$ has the law
    $\mathbf{P}^K(b(x),d(x),\alpha(x,x),N^K_x(0))$. Assume
    $N^K_x(0)\rightarrow n_x(0)$ in probability when
    $K\rightarrow+\infty$. Then, the sequence $(N^K_x)$ converges in
    probability on $[0,T]$ for the uniform norm to the deterministic
    function $n_x$ with initial condition $n_x(0)$ solution to
    \begin{equation}
      \label{eq:LV-monomorphic}
      \dot{n}_x=(b(x)-d(x)-\alpha(x,x)n_x)n_x.
    \end{equation}
  \item[\textmd{(b)}] Assume $\mu\equiv 0$ and
    $\nu^K_0=N^K_x(0)\delta_x+N^K_y(0)\delta_y$. Then, for any $t\geq
    0$, $\nu^K_t=N^K_x(t)\delta_x+N^K_y(t)\delta_y$, where
    $(N^K_x,N^K_y)$ has the law
    \begin{equation*}
      \mathbf{Q}^K(b(x),b(y),d(x),d(y),\alpha(x,x),\alpha(x,y),
      \alpha(y,x),\alpha(y,y),N^K_x(0),N^K_y(0)).
    \end{equation*}
    Assume $N^K_x(0)\rightarrow n_x(0)$ and $N^K_y(0)\rightarrow
    n_y(0)$ in probability when $K\rightarrow +\infty$. Then,
    $(N^K_x,N^K_y)$ converges in probability when $K\rightarrow
    +\infty$ on $[0,T]$ for the uniform norm to the deterministic
    function $(n_x,n_y)$ with initial condition $(n_x(0),n_y(0))$
    solution to
    \begin{equation}
      \label{eq:LV-IPS}
      \left\{\begin{array}{l}
          \dot{n}_x=(b(x)-d(x)
          -\alpha(x,x)n_x-\alpha(x,y)n_y)n_x \\
          \dot{n}_y=(b(y)-d(y)
          -\alpha(y,x)n_x-\alpha(y,y)n_y)n_y.
        \end{array}\right.
    \end{equation}
  \end{description}
\end{propbr}

Note that, under assumption~(A3), the logistic
equation~(\ref{eq:LV-monomorphic}) has two steady states, 0, unstable,
and $\bar{n}_x$, defined in~(\ref{eq:def-bar-n}), stable. The
system~(\ref{eq:LV-IPS}) has at least three steady states, $(0,0)$,
unstable, $(\bar{n}_x,0)$ and $(0,\bar{n}_y)$.


The assumption~(B) of Section~\ref{sec:models} is the mathematical
formulation of the impossibility of coexistence of two different
traits, in the sense that, starting in the neighborhood of the
equilibrium $(\bar{n}_x,0)$ of system~(\ref{eq:LV-IPS}), either its
solution converges to this equilibrium or to the equilibrium
$(0,\bar{n}_y)$. More precisely, the following proposition follows
fron an elementary analysis of system~(\ref{eq:LV-IPS}) (cf.\ e.g.\ 
Istas, 2000, pp.\ 25--27):
\begin{prop}
  \label{prop:LV}
  If $x$ and $y$ satisfy~(\ref{eq:hyp-B1}), then $(\bar{n}_x,0)$ is a
  stable steady state of~(\ref{eq:LV-IPS}). If $x$ and $y$
  satisfy~(\ref{eq:hyp-B2}), then $(\bar{n}_x,0)$ is an unstable
  steady state, $(0,\bar{n}_y)$ is stable, and any solution
  to~(\ref{eq:LV-IPS}) with initial state in $(\mathbb{R}_+^*)^2$
  converges to $(0,\bar{n}_y)$ when $t\rightarrow+\infty$.
\end{prop}

Let us now give the main ideas of the proof of
Theorem~\ref{thm:IPS-TSS}. It is based on two main ingredients: first,
when $\mu\equiv 0$ and $\nu^K_0$ is monomorphic with trait $x$, we
have seen in Proposition~\ref{prop:cv-mono-di}~(a) the convergence of
$\nu^K$ to $n(t)\delta_x$, where $n(t)$ is solution
to~(\ref{eq:LV-monomorphic}).  Any solution to this equation with
positive initial condition converges for large time to $\bar{n}_x$.
The large deviations estimates for this convergence will allow us to
show that the time during which the stochastic process stays in a
neighborhood of its limit (problem of exit from domain, Freidlin and
Wentzell, 1984) is of the order of $\exp(KV)$ with $V>0$. Now, when
$u_K$ is small, the process $\nu^K$ with a monomorphic initial
condition with trait $x$ is close to the same process with $\mu\equiv
0$, as long as no mutation occurs. Therefore, the left inequality
in~(\ref{eq:cond-mu-K}) will allow us to prove that, with high
probability, the first mutation event (occuring on the time scale
$t/Ku_K$) occurs before the total density drifts away from
$\bar{n}_x$.

The second ingredient of our proof is the study of the invasion of a
mutant trait $y$ that has just appeared in a monomorphic population
with trait $x$. This invasion can be divided in three steps
(Fig.~\ref{fig:inv-fix}), in a similar way as is done classically by
population geneticists dealing with selective sweeps (Kaplan et al.,
1989; Durrett and Schweinsberg, 2004):
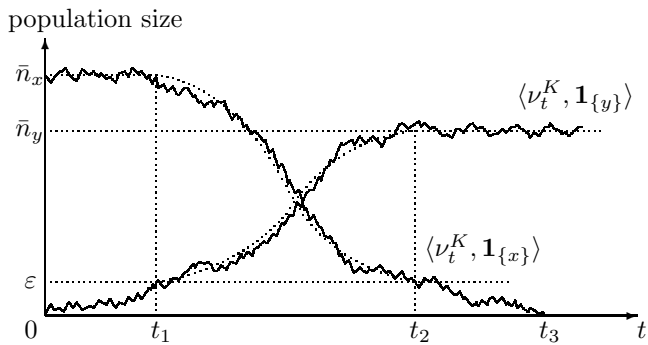
\begin{figure}[ht]
  \begin{center}
    \begin{picture}(350,180)(-20,-10)
      \put(0,0){\vector(1,0){320}} \put(0,0){\vector(0,1){150}}
      \put(-11,-12){0} \put(-11,15){$\varepsilon$}
      \dottedline{3}(0,18)(250,18)
      \put(-16,97){$\bar{n}_y$} \put(-16,127){$\bar{n}_x$}
      \dottedline{3}(0,100)(300,100) \dottedline{3}(0,130)(60,130)
      \put(-20,156){population size} \put(57,-12){$t_1$}
      \put(197,-12){$t_2$} \put(267,-12){$t_3$} \put(320,-12){$t$}
      \dottedline{3}(60,0)(60,130) \dottedline{3}(200,0)(200,100)
      \qbezier[27](60,130)(100,125)(130,70)
      \qbezier[25](130,70)(150,28)(200,18)
      \qbezier[25](60,18)(105,26)(130,55)
      \qbezier[25](130,55)(155,95)(200,100)
      \dottedline{0.5}(0,1)(3,6)(5,3)(8,4)(10,2)(12,6)(13,7)(15,4)(17,7)%
(20,3)(24,5)(26,4)(29,8)(33,6)(35,9)(38,4)(40,7)(43,6)(45,10)(47,8)(48,8)%
(50,12)(52,9)(55,13)(56,15)(57,14)(58,17)(59,15)(60,18)(61,19)(64,15)%
(66,16)(69,20)(71,20)(72,19)(75,23)(77,24)(80,28)(81,26)(83,29)(86,28)%
(88,29)(90,27)(91,24)(93,24)(95,27)(98,25)(100,29)(101,29)(103,28)(106,33)%
(107,31)(109,34)(110,34)(112,32)(113,36)(115,37)(118,40)(120,39)(123,44)%
(125,43)(126,46)(128,50)(129,49)(132,54)(134,54)(136,61)(138,59)(141,66)%
(143,68)(145,68)(148,75)(149,73)(151,79)(153,81)(154,84)(156,83)(158,89)%
(160,93)(161,92)(164,94)(166,95)(167,92)(170,97)(173,93)(175,95)(177,97)%
(181,93)(182,91)(184,96)(186,96)(189,98)(191,103)(193,104)(196,101)(199,103)%
(202,105)(203,102)(205,104)(208,101)(210,102)(213,99)(215,97)(217,98)%
(220,96)(222,101)(224,101)(226,99)(229,103)(231,104)(233,101)(236,102)%
(239,99)(240,97)(242,98)(244,98)(247,101)(249,99)(252,102)(254,102)(255,103)%
(257,101)(258,101)(261,98)(262,96)(264,99)(267,97)(270,102)(272,102)%
(273,104)(275,101)(278,100)(279,98)(281,99)(283,97)(288,102)(290,102)
      \put(255,115){$\langle\nu^K_t,\mathbf{1}_{\{y\}}\rangle$}
      \dottedline{0.5}(0,130)(3,126)(6,130)(8,132)(9,132)(11,134)(13,130)%
(15,128)(17,127)(20,133)(24,130)(26,132)(29,129)(30,131)(33,128)(35,129)%
(38,126)(40,131)(43,132)(45,129)(47,134)(48,133)(50,131)(52,132)(55,128)%
(56,130)(57,127)(58,125)(59,128)(60,129)(61,123)(64,125)(66,121)(69,124)%
(71,120)(72,119)(75,123)(77,123)(80,118)(81,119)(83,115)(84,116)(86,113)%
(88,115)(90,116)(91,113)(93,116)(95,115)(97,116)(98,112)(100,108)(101,106)%
(103,109)(106,106)(107,102)(109,103)(111,98)(113,100)(116,97)(118,94)%
(119,95)(121,88)(123,90)(124,86)(126,80)(128,82)(130,77)(131,71)(132,73)%
(135,66)(137,68)(139,62)(140,57)(141,58)(143,52)(145,54)(148,48)(149,48)%
(151,43)(153,42)(154,43)(156,39)(157,39)(159,32)(160,33)(162,28)(164,30)%
(166,25)(167,23)(170,26)(172,24)(174,23)(175,25)(177,24)(181,27)(182,26)%
(184,23)(186,24)(189,20)(191,22)(193,24)(194,21)(196,19)(198,20)(200,18)%
(202,15)(203,17)(205,20)(208,18)(210,21)(213,17)(215,14)(217,16)(220,12)%
(222,11)(224,13)(226,9)(227,10)(229,7)(231,5)(233,7)(234,5)(236,8)%
(239,6)(240,9)(242,10)(244,7)(247,10)(249,6)(250,7)(252,5)(254,8)%
(255,4)(257,6)(258,4)(261,2)(262,4)(264,3)(267,1)(269,2)(270,0)
      \put(205,32){$\langle\nu^K_t,\mathbf{1}_{\{x\}}\rangle$}
    \end{picture}
  \end{center}
  \caption{The three steps of the invasion of a
    mutant trait $y$ in a monomorphic population with trait $x$.}
  \label{fig:inv-fix}
\end{figure}
\begin{itemize}
\item Firstly, as long as the mutant population size
  $\langle\nu^K_t,\mathbf{1}_{\{y\}}\rangle$ (initially equal to
  $1/K$) is smaller than a fixed small $\varepsilon>0$ (before $t_1$
  in Fig.~\ref{fig:inv-fix}), the resident dynamics is very close to
  what it was before the mutation, so
  $\langle\nu^K_t,\mathbf{1}_{\{x\}}\rangle$ stays close to
  $\bar{n}_x$. Then, the death rate of a mutant individual is close to
  the constant $d(y)+\alpha(y,x)\bar{n}_x$. Since its birth rate is
  constant, equal to $b(y)$, we can approximate the mutant dynamics by
  a binary branching process. Therefore, the probability that
  $\langle\nu^K_t,\mathbf{1}_{\{y\}}\rangle$ reaches $\varepsilon$ is
  approximately equal to the probability that this branching process
  reaches $\varepsilon K$, which converges when $K\rightarrow+\infty$
  to its probability of non-extinction $[f(y,x)]_+/b(y)$.
\item Secondly, once $\langle\nu^K_t,\mathbf{1}_{\{y\}}\rangle$ has
  reached $\varepsilon$, by Proposition~\ref{prop:cv-mono-di}~(b), for
  large $K$, $\nu^K$ is close to the solution to~(\ref{eq:LV-IPS})
  with initial state $(\bar{n}_x,\varepsilon)$ (represented with
  dotted lines in Fig.~\ref{fig:inv-fix}) with high probability.  By
  Proposition~\ref{prop:LV}, this solution will be shown to reach the
  $\varepsilon$-neighborhood of $(0,\bar{n}_y)$ in finite time ($t_2$
  in Fig.~\ref{fig:inv-fix}).
\item Finally, once $\langle\nu^K_t,\mathbf{1}_{\{y\}}\rangle$ is
  close to $\bar{n}_y$ and $\langle\nu^K_t,\mathbf{1}_{\{x\}}\rangle$
  is small, $K\langle\nu^K_t,\mathbf{1}_{\{x\}}\rangle$ can be
  approximated, in a similar way as in the first step, by a binary
  branching process, which is subcritical and hence gets extinct
  a.s.~in finite time ($t_3$ in Fig.~\ref{fig:inv-fix}).
\end{itemize}
We will see in Sections~\ref{sec:exit-dom} and~\ref{sec:br-proc} that
the time needed to complete the first and third steps is proportional
to $\log K$, whereas the time needed for the second step is bounded.
Therefore, since the time between two mutations is of the order of
$1/Ku_K$, the right inequality in~(\ref{eq:cond-mu-K}) will allow us
to prove that, with high probability, the three steps above are
completed before a new mutation occurs.

\begin{rem}
  \label{rem:init-cond}
  As observed by Metz et al.\ (1996), the biologial heuristics leading
  to the TSS model extend to the case of polymorphic initial
  condition, where the population is composed of a finite number of
  distinct traits (see also Champagnat, 2004). Our mathematical method
  can also be extended easily to $n$-morphic initial conditions,
  except for one difficulty: one has to replace assumption~(B) by
  another assumption stating that, for any $n$, any solution to the
  $n$-morphic logistic systems generalizing~(\ref{eq:LV-IPS})
  converges to an equilibrium (as in Proposition~\ref{prop:LV}), and
  that the equilibria of these systems are non-degenerate, in the
  sense that the branching processes in the first and third steps
  above are not critical, or, equivalently, that a first-order linear
  analysis of these equilibria allows to determine their stability.
  Then, one could construct a polymorphic TSS model in which the
  number of coexisting traits is not fixed.  However, the asymptotic
  analysis of $n$-dimensional logistic systems is non-trivial and may
  exhibit cycles or chaos, except when $n=1$ or $2$, and analytical
  assumptions ensuring the condition above are difficult to find.
\end{rem}

Section~\ref{sec:birth-death-proc} will provide the large deviations
and branching process results needed to make formal the previous
heuristics. We will also prove several comparison results between
$\langle\nu^K_t,\mathbf{1}\rangle$ and the birth and death processes
of Definition~\ref{def:birth-death-proc}. In Section~\ref{sec:proof},
the proof of Theorem~\ref{thm:IPS-TSS} is achieved by computing, for
any $t$, the limit law of $\nu^K_{t/Ku_K}$ according to the random
number of mutations having occured between 0 and $t/Ku_K$.

\medskip

\noindent\textbf{Notations}
\begin{itemize}
\item $\lceil a\rceil$ denotes the first integer greater or equal to
  $a$, and $\lfloor a\rfloor$ denotes the integer part of $a$.
\item For any $K\geq 1$ and $\nu\in{\cal M}^K$, we will denote by
  $\mathbf{P}^K_{\nu}$ the law of the process $\nu^K$ generated
  by~(\ref{eq:generator-renormalized-IPS}) with initial state $\nu$,
  and by $\mathbf{E}^K_{\nu}$ the expectation with respect to
  $\mathbf{P}^K_{\nu}$.
\item The convergence in probability of finite dimensional random
  variables will be denoted by $\overset{\cal P}{\rightarrow}$.
\item We will denote by ${\cal L}(Z)$ the law of the stochastic
  process $(Z_t,t\geq 0)$.
\item We will denote by $\preceq$ the following stochastic domination
  relation: if $\mathbf{Q}_1$ and $\mathbf{Q}_2$ are the laws of
  $\mathbb{R}$-valued processes, we will write
  $\mathbf{Q}_1\preceq\mathbf{Q}_2$ if we can construct on the same
  probability space $(\Omega,{\cal F},\mathbf{P})$ two processes $X^1$
  and $X^2$ such that ${\cal L}(X^i)=\mathbf{Q}_i$ ($i=1,2$) and
  $\forall t\geq 0$, $\forall\omega\in\Omega$, $X^1_t(\omega)\leq
  X^2_t(\omega)$.
\item Finally, if $X^1$ and $X^2$ are two random processes and $T$ is
  a random time constructed on the same probability space as $X^1$, we
  will write $X^1_t\preceq X^2_t$ for $t\leq T$ (resp.\ $X^2_t\preceq
  X^1_t$ for $t\leq T$) if we can construct a process $\hat{X}^2$ on
  the same probability space as $X^1$, such that ${\cal
    L}(\hat{X}^2)={\cal L}(X^2)$ and $\forall t\leq T$,
  $\forall\omega\in\Omega$, $X^1_t(\omega)\leq\hat{X}^2_t(\omega)$
  (resp.\ $\hat{X}^2_t(\omega)\leq X^1_t(\omega)$).
\end{itemize}

\section{Birth and death processes}
\label{sec:birth-death-proc}

We will collect in this section various results on the birth and
death processes that appeared in
Definition~\ref{def:birth-death-proc}.

\subsection{Comparison results}
\label{sec:compar}

The following theorem gives various stochastic domination results.
\begin{thmbr}
  \label{thm:compar}
  \begin{description}
  \item[\textmd{(a)}] Assume~(A). For any $K\geq 1$ and any
    $\mathbb{L}^1$ initial condition $\nu^K_0$ of the process $\nu^K$,
    \begin{equation*}
      {\cal L}(\langle\nu^K,\mathbf{1}\rangle)\preceq
      \mathbf{P}^K(2\bar{b},0,\underline{\alpha},
      \langle\nu^K_0,\mathbf{1}\rangle).
    \end{equation*}
  \item[\textmd{(b)}] With the same assumptions as in~(a), let
    $A^K_t$ denote the number of mutations occuring in $\nu^K$ between
    times 0 and $t$, and let $a,a_1,a_2\geq 0$. Then, for
    $t\leq\inf\{s\geq 0:\langle\nu^K_s,\mathbf{1}\rangle\geq a\}$,
    \begin{equation*}
      A^K_t\preceq B^K_t,
    \end{equation*}
    where $B^K$ is a Poisson process with parameter $Ku_Ka\bar{b}$.
    
    If moreover $\nu^K_0=\langle\nu^K_0,\mathbf{1}\rangle\delta_x$,
    define $\tau_1=\inf\{t\geq 0: A^K_t=1\}$ (the first mutation
    time). Then, for $t\leq\tau_1\wedge\inf\{s\geq
    0:\langle\nu^K_s,\mathbf{1}\rangle\not\in[a_1,a_2]\}$,
    \begin{equation}
      \label{eq:pf-compar}
      B^K_t\preceq A^K_t\preceq C^K_t,
    \end{equation}
    where $B^K$ and $C^K$ are Poisson processes with respective
    parameter $Ku_Ka_1\mu(x)b(x)$ and $Ku_Ka_2\mu(x)b(x)$.
  \item[\textmd{(c)}] Fix $K\geq 1$ and take $b,d,\alpha,z$ as in
    Definition~\ref{def:birth-death-proc}~(a). Then, for any
    $\varepsilon_1,\varepsilon_2,\varepsilon_3\geq 0$ and any
    $\mathbb{N}/K$-valued random variable $\varepsilon_4$,
    \begin{equation*}
      \mathbf{P}^K(b,d+\varepsilon_2,\alpha+\varepsilon_3,z)\preceq
      \mathbf{P}^K(b+\varepsilon_1,d,\alpha,z+\varepsilon_4).
    \end{equation*}
  \item[\textmd{(d)}] Let $(Z^1,Z^2)$ be a stochastic process with law
    \begin{equation*}
      \mathbf{Q}^K(b_1,b_2,d_1,d_2,\alpha_{11},\alpha_{12},
      \alpha_{21},\alpha_{22},z_1,z_2)
    \end{equation*}
    where the parameters are as in
    Definition~\ref{def:birth-death-proc}~(b). Fix $a>0$ and define
    $T=\inf\{t\geq 0,Z^2\geq a\}$. Then, for $t\leq T$,
    \begin{gather*}
      M^1_t\preceq Z^1_t\preceq M^2_t, \\
      \mbox{where}\quad {\cal
        L}(M^1)=\mathbf{P}^K(b_1,d_1+a\alpha_{12},\alpha_{11},z_1) \\
      \mbox{and}\quad {\cal L}(M^2)=\mathbf{P}^K(b_1,d_1,\alpha_{11},z_1).
    \end{gather*}
  \item[\textmd{(e)}] Take $(Z^1,Z^2)$ as above, fix $0\leq a_1<a_2$
    and $a>0$, and define $T=\inf\{t\geq 0,Z^1\not\in[a_1,a_2]\mbox{\ 
      or\ }Z^2\geq a\}$. Then, for $t\leq T$,
    \begin{gather*}
      M^1_t\preceq Z^2_t\preceq M^2_t, \\
      \mbox{where}\quad {\cal
        L}(M^1)=\mathbf{P}^K(b_2,d_2+a_2\alpha_{21}+a\alpha_{22},0,z_2) \\
      \mbox{and}\quad {\cal
        L}(M^2)=\mathbf{P}^K(b_2,d_2+a_1\alpha_{21},0,z_2).
    \end{gather*}
  \end{description}
\end{thmbr}

\begin{rem}
  \label{rem:compar}
  Point~(a) explains why it is necessary to combine simultaneously the
  limits $K\rightarrow +\infty$ and $u_K\rightarrow 0$ in order to
  obtain the TSS process in Theorem~\ref{thm:IPS-TSS}. The limit
  $K\rightarrow +\infty$ taken alone leads to a deterministic dynamics
  (Fournier and M\'el\'eard, 2003), so making the rare mutations limit
  afterwards cannot lead to a stochastic process.  Conversely, taking
  the limit of rare mutations without making the population larger
  would lead to an immediate extinction of the population in the
  mutations time scale, because the stochastic domination of
  Theorem~\ref{thm:compar}~(a) is independent of $u_K$ and
  $\mu(\cdot)$, and because a process $Z$ with law
  $\mathbf{P}^K(2\bar{b},0,\underline{\alpha},\gamma_K/K)$ gets a.s.\ 
  extinct in finite.
\end{rem}

Before proving Theorem~\ref{thm:compar}, let us deduce from Point~(a)
the Lemma~\ref{lem:moment} stated in Section~\ref{sec:models}.

\paragraph{Proof of Lemma~\ref{lem:moment}}
By Theorem~\ref{thm:compar}~(a), it suffices to prove that
\begin{equation*}
  \sup_{K\geq 1}\:\sup_{t\geq 0}\mathbf{E}((Z^K_t)^p)<+\infty,
\end{equation*}
where ${\cal L}(Z^K)=\mathbf{P}^K(2\bar{b},0,\underline{\alpha},
z^K_0)$ when $\:\sup_{K\geq 1} \mathbf{E}((z^K_0)^p)<+\infty$.

Let us define $v^k_t=\mathbf{P}(Z^K_t=k/K)$. Then
\begin{align*}
  \frac{d}{dt}\mathbf{E}((Z^K_t)^p) & =\sum_{k\geq
    1}\left(\frac{k}{K}\right)^p\frac{dv_t^k}{dt} \\
  & \ =\frac{1}{K^p}\sum_{k\geq 1}k^p\left[2\bar{b}(k-1)v^{k-1}_t
    +\underline{\alpha}\frac{(k+1)^2}{K}v^{k+1}_t
    -k\left(2\bar{b}+\underline{\alpha}\frac{k}{K}\right)v^k_t\right] \\
  & \ =\frac{1}{K^p}\sum_{k\geq
    1}\left[2\bar{b}\left(\left(1+\frac{1}{k}\right)^p-1\right)
    +\underline{\alpha}\frac{k}{K}
    \left(\left(1-\frac{1}{k}\right)^p-1\right)\right]k^{p+1}v^k_t.
\end{align*}
Now, for $k/K>4\bar{b}/\underline{\alpha}$, the quantity inside the
square brackets in the last expression can be upper bounded by
$-2\bar{b}[3-2(1-1/k)^p-(1+1/k)^p]$, which is equivalent to
$-2\bar{b}p/k$ when $k\rightarrow+\infty$. Therefore, there exists a
constant $k_0$ that can be assumed bigger than
$4\bar{b}/\underline{\alpha}$ such that, for any $k\geq k_0$,
$-2\bar{b}[3-2(1-1/k)^p-(1+1/k)^p]\leq -\bar{b}p/k$. Then, using the
fact that $(1+x)^p-1\leq x(2^p-1)$ for any $x\in[0,1]$,
we can write
\begin{align*}
  \frac{d}{dt}\mathbf{E}((Z^K_t)^p) & \leq\sum_{k=1}^{Kk_0-1}
  2\bar{b}(2^p-1)\left(\frac{k}{K}\right)^pv^k_t
  -\sum_{k\geq Kk_0}\bar{b}p\left(\frac{k}{K}\right)^pv^k_t \\
  & \leq 2\bar{b}(2^p-1)k_0^p+\bar{b}pk_0^p
  -\bar{b}p \mathbf{E}((Z^K_t)^p).
\end{align*}
Writing $C=(2(2^p-1)+p)k_0^p/p$, this differential inequality
solves as
\begin{equation*}
  \mathbf{E}((Z^K_t)^p)\leq C+[\mathbf{E}((z_0^K)^p)-C]e^{-\bar{b}pt},
\end{equation*}
which gives the required uniform bound.\hfill$\Box$

\paragraph{Proof of Theorem~\ref{thm:compar}}
The proof is essentially intuitive if one computes upper and lower
bounds of the birth and death rates for each processes considered in
the statement of the theorem. We will simply give the explicit
construction of the process $\nu^K$, and the proof
of~(\ref{eq:pf-compar}) as an example. We leave the remaining
comparison results to the reader.

We will use the construction of the process $\nu^K$ given by Fournier
and M\'el\'eard (2003): let $(\Omega,{\cal F},\mathbf{P})$ be a
sufficiently large probability space, and consider on this space the
following five independent random objects:
\begin{description}
\item[\textmd{(i)}] a ${\cal M}^K$-valued random variable $\nu^K_0$
  (the initial distribution),
\item[\textmd{(ii)}] a Poisson point measure $N_1(ds,di,dv)$ on
  $[0,\infty[\times\mathbb{N}\times[0,1]$ with intensity measure
  $q_1(ds,di,dv)=\bar{b}\:ds\sum_{k\geq 1}\delta_k(di)dv$ (the birth
  without mutation Poisson point measure),
\item[\textmd{(iii)}] a Poisson point measure $N_2(ds,di,dh,dv)$ on
  $[0,\infty[\times\mathbb{N}\times\mathbb{R}^l\times[0,1]$ with
  intensity measure $q_2(ds,di,dh,dv)=\bar{b}\:ds\sum_{k\geq
    1}\delta_k(di)\bar{m}(h)dhdv$ (the birth with mutation Poisson point
  measure),
\item[\textmd{(iv)}] a Poisson point measure $N_3(ds,di,dv)$ on
  $[0,\infty[\times\mathbb{N}\times[0,1]$ with intensity measure
  $q_3(ds,di,dv)=\bar{d}\:ds\sum_{k\geq 1}\delta_k(di)dv$ (the natural
  death Poisson point measure),
\item[\textmd{(v)}] a Poisson point measure $N_4(ds,di,dj,dv)$ on
  $[0,\infty[\times\mathbb{N}\times\mathbb{N}\times[0,1]$ with
  intensity measure $q_4(ds,di,dj,dv)=(\bar{\alpha}/K)ds\sum_{k\geq
    1}\delta_k(di)\sum_{m\geq 1}\delta_m(dj)dv$ (the competition death
  Poisson point measure).
\end{description}
We will also need the following function, solving the purely
notational problem of associating a number to each individual in the
population: for any $K\geq 1$, let $H=(H^1,\ldots,H^k,\ldots)$ be the
map from ${\cal M}^K$ into $(\mathbb{R}^l)^{\mathbb{N}}$ defined by
\begin{equation*}
  H\left(\frac{1}{K}\sum_{i=1}^n\delta_{x_i}\right)=(x_{\sigma(1)},\ldots,
      x_{\sigma(n)},0,\ldots,0,\ldots),
\end{equation*}
where $x_{\sigma(1)}\curlyeqprec\ldots\curlyeqprec x_{\sigma(n)}$ for
the lexicographic order $\curlyeqprec$ on $\mathbb{R}^l$. For
convenience, we have omitted in our notation the dependence of $H$ and
$H^i$ on $K$.

Then a process $\nu^K$ with generator $L^K$ and initial state
$\nu_0^K$ can be constructed as follows: for any $t\geq 0$,
\begin{align}
  \nu^K_t & =\nu^K_0+\int_0^t\int_{\mathbb{N}}\int_0^1\mathbf{1}_{\{i\leq
    K\langle\nu^K_{s-},\mathbf{1}\rangle\}}
  \frac{\delta_{H^i(\nu^K_{s-})}}{K} \notag \\
  & \qquad\qquad
  \mathbf{1}_{\left\{v\leq\frac{[1-u_K\mu(H^i(\nu^K_{s-}))]
        b(H^i(\nu^K_{s-}))}{\bar{b}}\right\}}N_1(ds,di,dv) \notag \\
  & +\int_0^t\int_{\mathbb{N}}\int_{\mathbb{R}^l}\int_0^1
  \mathbf{1}_{\{i\leq K\langle\nu^K_{s-},\mathbf{1}\rangle\}}
  \frac{\delta_{H^i(\nu^K_{s-})+h}}{K} \notag \\
  & \qquad\qquad
  \mathbf{1}_{\left\{v\leq\frac{u_K\mu(H^i(\nu^K_{s-}))
        b(H^i(\nu^K_{s-}))}{\bar{b}}\frac{m(H^i(\nu^K_{s-}),h)}{\bar{m}(h)}
    \right\}}N_2(ds,di,dh,dv) \notag \\
  & -\int_0^t\int_{\mathbb{N}}\int_0^1\mathbf{1}_{\{i\leq
    K\langle\nu^K_{s-},\mathbf{1}\rangle\}}\frac{\delta_{H^i(\nu^K_{s-})}}{K}
  \mathbf{1}_{\left\{v\leq\frac{d(H^i(\nu^K_{s-}))}{\bar{d}}
    \right\}}N_3(ds,di,dv) \notag \\
  & -\int_0^t\int_{\mathbb{N}}\int_{\mathbb{N}}\int_0^1
  \mathbf{1}_{\{i\leq K\langle\nu^K_{s-},\mathbf{1}\rangle\}}
  \mathbf{1}_{\{j\leq K\langle\nu^K_{s-},\mathbf{1}\rangle\}}
  \frac{\delta_{H^i(\nu^K_{s-})}}{K} \notag \\
  & \qquad\qquad
  \mathbf{1}_{\left\{v\leq\frac{\alpha(H^i(\nu^K_{s-}),H^j(\nu^K_{s-}))}
      {\bar{\alpha}}\right\}}N_4(ds,di,dj,dv). \label{eq:def-XK}
\end{align}
Although this formula is quite complicated, the principle is
simple: for each type of event, the corresponding Poisson point
process jumps faster than $\nu^K$ has to. We decide whether a jump of
the process $\nu^K$ occurs by comparing $v$ to a quantity related to
the rates of the various events. The indicator functions involving $i$
and $j$ ensures that the $i^{\mbox{\footnotesize{th}}}$ and
$j^{\mbox{\footnotesize{th}}}$ individuals are alive in the population
(because $K\langle\nu^K_t,\mathbf{1}\rangle$ is the number of
individuals in the population at time $t$).

Under~(A1), (A2) and the assumption that
$\mathbf{E}(\langle\nu^K_0,\mathbf{1}\rangle)<\infty$, Fournier and
M\'el\'eard (2003) prove the existence and uniqueness of the solution
to~(\ref{eq:def-XK}), and that this solution is a Markov process with
infinitesimal generator~(\ref{eq:generator-renormalized-IPS}).

Now, let us come to the proof of~(\ref{eq:pf-compar}). The process
$A^K$ can be written as
\begin{multline*}
  A^K_t:=\int_0^t\int_{\mathbb{N}} \int_{\mathbb{R}^l}\int_0^1
  \mathbf{1}_{\{i\leq K\langle\nu^K_{s-},\mathbf{1}\rangle\}}\times \\
  \times\mathbf{1}_{\left\{v\leq\frac{u_K\mu(H^i(\nu^K_{s-}))
        b(H^i(\nu^K_{s-}))}{\bar{b}}\frac{m(H^i(\nu^K_{s-}),h)}{\bar{m}(h)}
    \right\}}N_2(ds,di,dh,dv).
\end{multline*}
In the case where $\nu^K_0=\langle\nu^K_0,\mathbf{1}\rangle\delta_x$,
as long as $t<\tau_1$,
$\nu^K_t=\langle\nu^K_t,\mathbf{1}\rangle\delta_x$. Therefore, for
$t\leq\tau_1\wedge\inf\{s\geq
0:\langle\nu^K_s,\mathbf{1}\rangle\not\in[a_1,a_2]\}$,
\begin{multline}
  \label{eq:pf-compar-2}
  \int_0^t\int_{\mathbb{N}}\int_{\mathbb{R}^l}\int_0^1
  \mathbf{1}_{\{i\leq Ka_1\}}
  \mathbf{1}_{\left\{v\leq\frac{u_K\mu(x)b(x)}{\bar{b}}
      \frac{m(x,h)}{\bar{m}(h)}\right\}}N_2(ds,di,dh,dv)\leq A^K_t \\ \leq
  \int_0^t\int_{\mathbb{N}}\int_{\mathbb{R}^l}\int_0^1
  \mathbf{1}_{\{i\leq Ka_2\}}
  \mathbf{1}_{\left\{v\leq\frac{u_K\mu(x)b(x)}{\bar{b}}
      \frac{m(x,h)}{\bar{m}(h)}\right\}}N_2(ds,di,dh,dv).
\end{multline}
Since the intensity measure of $N_2$ is
\begin{equation*}
  q_2(ds,di,dh,dv)=\bar{b}\:ds\sum_{k\geq 1}\delta_k(di)\bar{m}(h)dhdv,  
\end{equation*}
the left-hand side and the right-hand side of~(\ref{eq:pf-compar-2})
are Poisson processes with parameters $Ku_Ka_1\mu(x)b(x)$ and
$Ku_Ka_2\mu(x)b(x)$, respectively.\hfill$\Box$

\subsection{Problem of exit from a domain}
\label{sec:exit-dom}

Let us give some results on $\mathbf{P}^K(b,d,\alpha,z)$ when $\alpha>0$.
Points~(a) and~(b) of the following theorem strengthen
Proposition~\ref{prop:cv-mono-di}, and point~(c) studies the problem
of exit from a domain.
\begin{thmbr}
  \label{thm:exit-time}
  \begin{description}
  \item[\textmd{(a)}] Let $\alpha,T>0$ and $b,d\geq 0$, let $C$ be a
    compact subset of $\mathbb{R}_+^*$, and write
    $\mathbf{P}^K_z=\mathbf{P}^K(b,d,\alpha,z)$ for $z\in\mathbb{N}/K$. Let
    $\phi_z$ denote the solution to
    \begin{equation}
      \label{eq:b-d-proc-lim-mono}
      \dot{\phi}=(b-d-\alpha\phi)\phi
    \end{equation}
    with initial condition $\phi_z(0)=z$. Then
    \begin{equation*}
      r:=\inf_{z\in C}\ \inf_{0\leq t\leq T}|\phi_z(t)|>0\mbox{\ and\ }
      R:=\sup_{z\in C}\ \sup_{0\leq t\leq T}|\phi_z(t)|<+\infty.
    \end{equation*}
    Moreover, for any $\delta<r$,
    \begin{equation}
      \label{eq:b-d-proc-LD}
      \lim_{K\rightarrow +\infty}\ \sup_{z\in C}
      \mathbf{P}^K_z\biggl(\sup_{0\leq t\leq T}|w_t-\phi_z(t)|
      \geq\delta\biggr)=0,
    \end{equation}
    where $w_t$ is the canonical process on
    $\mathbb{D}(\mathbb{R}_+,\mathbb{R})$.
  \item[\textmd{(b)}] Let $T,\alpha_{ij}>0$ and $b_i,d_i\geq 0$
    ($i,j\in\{1,2\}$), let $C$ be a compact subset of
    $(\mathbb{R}_+^*)^2$, and write $\mathbf{Q}^K_{z_1,z_2}=
    \mathbf{Q}^K(b_1,b_2,d_1,d_2,\alpha_{11},\alpha_{12},
    \alpha_{21},\alpha_{22},z_1,z_2)$
    for $z_1$ and $z_2$ in $\mathbb{N}/K$. Let
    $\phi_{z_1,z_2}=(\phi^1_{z_1,z_2},\phi^2_{z_1,z_2})$ denote the
    solution to
    \begin{equation*}
      \left\{
        \begin{array}{l}
          \dot{\phi}^1 = (b_1-d_1-\alpha_{11}\phi^1-\alpha_{12}\phi^2)\phi^1 \\
          \dot{\phi}^2 = (b_2-d_2-\alpha_{21}\phi^1-\alpha_{22}\phi^2)\phi^2
        \end{array}
      \right.
    \end{equation*}
    with initial conditions $\phi^1_{z_1,z_2}(0)=z_1$ and
    $\phi^2_{z_1,z_2}(0)=z_2$. Then
    \begin{equation}
      \label{eq:def-r}
      r:=\inf_{z\in C}\ \inf_{0\leq t\leq
      T}\|\phi_{z_1,z_2}(t)\|>0\mbox{\ and\ }
      \sup_{z\in C}\ \sup_{0\leq t\leq T}\|\phi_{z_1,z_2}(t)\|<+\infty.
    \end{equation}
    Moreover, for any $\delta<r$,
    \begin{equation*}
      \lim_{K\rightarrow +\infty}\sup_{z\in C}
      \mathbf{Q}^K_{z_1,z_2}(\sup_{0\leq t\leq T}
      \|\hat{w}_t-\phi_{z_1,z_2}(t)\|\geq\delta)=0,
    \end{equation*}
    where $\hat{w}_t=(\hat{w}^1_t,\hat{w}^2_t)$ is the canonical process on
    $\mathbb{D}(\mathbb{R}_+,\mathbb{R}^2)$.
  \item[\textmd{(c)}] Let $b,\alpha>0$ and $0\leq d<b$. Observe that
    $(b-d)/\alpha$ is the unique stable steady state
    of~(\ref{eq:b-d-proc-lim-mono}).
    Fix $0<\eta_1<(b-d)/\alpha$ and $\eta_2>0$, and define on
    $\mathbb{D}(\mathbb{R}_+,\mathbb{R})$
    \begin{equation*}
      T^K=\inf\left\{t\geq
        0:w_t\not\in\left[\frac{b-d}{\alpha}-\eta_1,
          \frac{b-d}{\alpha}+\eta_2\right]\right\}. 
    \end{equation*}
    Then, there exists $V>0$ such that, for any compact subset $C$ of
    $](b-d)/\alpha-\eta_1,(b-d)/\alpha+\eta_2[$,
    \begin{equation}
      \label{eq:exit-time}
      \lim_{K\rightarrow+\infty}\ \sup_{z\in C}
      \mathbf{P}^K_z(T^K<e^{KV})=0.
    \end{equation}
  \end{description}
\end{thmbr}

\paragraph{Proof of~(a) and~(b)}
Observe that any solution to~(\ref{eq:b-d-proc-lim-mono}) with
positive initial condition is bounded ($\dot{\phi}<0$ as soon as
$\phi>(b-d)/\alpha$). This implies that $R<\infty$.  Moreover, a
solution to~(\ref{eq:b-d-proc-lim-mono}) can be written as
\begin{equation*}
  \phi(t)=\phi(0)\exp\left(\int_0^t(b-d-\alpha\phi(s))ds\right)
  \geq\phi(0)\exp((b-d-\alpha R)t),
\end{equation*}
which implies that $r>0$.

Equation~(\ref{eq:b-d-proc-LD}) is a consequence of large deviations
estimates for the sequence of laws $(\mathbf{P}^K_z)_{K\geq 1}$. As
can be seen in Theorem~10.2.6 in Chap.~10 of Dupuis and Ellis (1997),
a large deviations principle on $[0,T]$ with a good rate function
$I_T$ holds for $\mathbb{Z}/K$-valued Markov jump processes with
transition rates
\begin{equation*}
  \begin{array}{ll}
    Kp(i/K) & \mbox{from\ }i/K\mbox{\ to\ }(i+1)/K, \\
    Kq(i/K) & \mbox{from\ }i/K\mbox{\ to\ }(i-1)/K,
  \end{array}
\end{equation*}
where $p$ and $q$ are functions defined on $\mathbb{R}$ and with
positive values, bounded, Lipschitz and uniformly bounded away from 0.
The rate function $I_T$ writes
\begin{equation}
  \label{eq:def-I-Dupuis-Ellis}
  I_T(\phi)=\left\{
    \begin{array}{ll}
      \displaystyle{\int_0^TL(\phi(t),\dot{\phi}(t))dt} & \mbox{if $\phi$
          is absol.\ cont.\ on $[0,T]$} \\
        +\infty & \mbox{otherwise}
    \end{array}\right.
\end{equation}
for some function $L:\mathbb{R}^2\rightarrow\mathbb{R}_+$ such that
$L(y,z)=0$ if and only if $z=p(y)-q(y)$. Therefore, $I_T(\phi)=0$ if
and only if $\phi$ is absolutely continuous and
\begin{equation}
  \label{eq:lim-det-p-q}
  \dot{\phi}=p(\phi)-q(\phi).
\end{equation}

Moreover, this large deviation is uniform with respect to the
initial condition. This means that, if $\mathbf{R}^K_z$ denotes the
law of this process with initial condition $z$, for any compact set
$C\subset\mathbb{R}$, for any closed set $F$ and any open set $G$ of
$\mathbb{D}([0,T],\mathbb{R})$,
\begin{align}
  \liminf_{K\rightarrow+\infty}\frac{1}{K}\log \inf_{z\in
    C}\mathbf{R}^K_z(G) & \geq -\sup_{z\in C}\ \inf_{\psi\in G,\ 
    \psi(0)=z}I_T(\psi) \label{eq:LDP-LB} \\
  \mbox{and}\quad\limsup_{K\rightarrow+\infty}\frac{1}{K}\log\sup_{z\in
    C}\mathbf{R}^K_z(F) & \leq-\inf_{\psi\in F,\ \psi(0)\in
    C}I_T(\psi). \label{eq:LDP-UB}
\end{align}

Our birth and death process does not satisfy these
asumptions. However, if we define
\begin{gather*}
  p(z)=b\chi(z) \quad\mbox{and}\quad q(z)=d\chi(z)+\alpha\chi(z)^2, \\
  \mbox{where}\quad \chi(z)=z \mbox{\ if\ }z\in [r-\delta,R+\delta];\ 
  r-\delta \mbox{\ if\ }z<r-\delta;\ R+\delta\mbox{\ if\ }z>R+\delta,
\end{gather*}
then $\mathbf{R}^K_z=\mathbf{P}^K_z$ on the time interval $[0,\tau]$,
where $\tau=\inf\{t\geq 0,w_t\not\in[r-\delta,R+\delta]\}$, and $p$
and $q$ satisfy the assumptions above. Therefore,
by~(\ref{eq:LDP-UB}),
\begin{gather*}
  \limsup_{K\rightarrow +\infty}\frac{1}{K}\log\sup_{z\in C}
  \mathbf{P}^K_z\biggl(\sup_{0\leq t\leq T}|w_t-\phi_z(t)|
  \geq\delta\biggr)\leq -\inf_{\psi\in F^{\delta}}I_T(\psi),
  \quad\mbox{where} \\
  F^{\delta}:=
  \bigl\{\psi\in\mathbb{D}([0,T],\mathbb{R}):\psi(0)\in C\mbox{\ and
    $\exists t\in[0,T]$,\ }|\psi(t)-\phi_{\psi(0)}(t)|\geq\delta\bigr\}
\end{gather*}
By the continuity of the flow of~(\ref{eq:lim-det-p-q}) (which is a
classical consequence of the fact that $z\mapsto p(z)-q(z)$ is
Lipschitz and of Gronwall's Lemma), the set $F^{\delta}$ is closed.
Since $I_T$ is a good rate function, the infimum of $I_T$ over this
set is attained at some function belonging to $F^{\delta}$, which
cannot be a solution to~(\ref{eq:lim-det-p-q}), and thus is non-zero.
This ends the proof of~(\ref{eq:b-d-proc-LD}).

The proof of~(b) can be made in a very similar way.\hfill$\Box$

\paragraph{Proof of~(c)}
Define the function $\chi$ on $\mathbb{R}$ by $\chi(z)=z$ if
$z\in[(b-d)/\alpha-\eta_1,(b-d)/\alpha+\eta_2]$,
$\chi(z)=(b-d)/\alpha-\eta_1$ for $z<(b-d)/\alpha-\eta_1$ and
$\chi(z)=(b-d)/\alpha+\eta_2$ for $z>(b-d)/\alpha-\eta_2$. As in the
proof of~(a), we can construct from the functions $p(z)=b\chi(z)$ and
$q(z)=d\chi(z)+\alpha\chi(z)^2$ a family of laws $(\mathbf{R}^K_z)$ such
that $\mathbf{R}^K_z=\mathbf{P}^K_z$ on the time interval $[0,T^K]$,
and such that~(\ref{eq:LDP-LB}) and~(\ref{eq:LDP-UB}) hold for the
good rate function $I_T$ defined in~(\ref{eq:def-I-Dupuis-Ellis}).

Observe that any solution to~(\ref{eq:lim-det-p-q}) are monotonous and
converge to $(b-d)/\alpha$ when $t\rightarrow +\infty$. Therefore, the
following estimates for the time of exit from an attracting domain are
classical (Freidlin and Wentzell, 1984, Chap.~5, Section~4): there
exists $\bar{V}\geq 0$ such that, for any $\delta>0$,
\begin{equation*}
  \lim_{K\rightarrow+\infty}\inf_{z\in C}
  \mathbf{R}^K_z\left(e^{K(\bar{V}-\delta)}
    <T^K<e^{K(\bar{V}+\delta)}\right)=1,
\end{equation*}
which implies~(\ref{eq:exit-time}) if we can prove that $\bar{V}>0$.

The constant $\bar{V}$ is obtained as follows (see Freidlin and
Wentzell, 1984, pp.\ 108--109): for any $y,z\in\mathbb{R}$, define
\begin{equation*}
  V(y,z):=\inf_{t>0,\ \varphi(0)=y,\ \varphi(t)=z}I_t(\varphi).
\end{equation*}
Then
\begin{equation*}
  \bar{V}:=V\left(\frac{b-d}{\alpha},\frac{b-d}{\alpha}-\eta_1\right)\wedge
  V\left(\frac{b-d}{\alpha},\frac{b-d}{\alpha}+\eta_2\right).
\end{equation*}

Now, Theorem~5.4.3.\ of Freidlin and Wentzell (1984) states that, for
any $y,z\in\mathbb{R}$, the infimum defining $V(y,z)$ is attained at
some function $\phi$ linking $y$ to $z$, in the sense that, either
there exists an absolutely continuous function $\phi$ defined on
$[0,T]$ for some $T>0$ such that $\phi(0)=y$, $\phi(T)=z$ and
$V(y,z)=I_T(\phi)=\int_0^TL(\phi(t),\dot{\phi}(t))dt$, or there exists
an absolutely continuous function $\phi$ defined on $]-\infty,T]$ for
some $T>-\infty$ such that $\lim_{t\rightarrow -\infty}\phi(t)=y$,
$\phi(T)=z$ and $V(y,z)=\int_{-\infty}^TL(\phi(t),\dot{\phi}(t))dt$.

Since any solution to~(\ref{eq:lim-det-p-q}) is decreasing as long as
it stays in $[(b-d)/\alpha,+\infty[$, a function $\phi$ defined on
$[0,T]$ or $]-\infty,T]$ linking $(b-d)/\alpha$ to
$(b-d)/\alpha+\eta_2$ cannot be a solution to~(\ref{eq:lim-det-p-q}),
and thus $V((b-d)/\alpha,(b-d)/\alpha+\eta_2)>0$. Similarly,
$V((b-d)/\alpha,(b-d)/\alpha-\eta_1)>0$, and so $\bar{V}>0$, which
concludes the proof of Theorem~\ref{thm:exit-time}.\hfill$\Box$

\subsection{Some results on branching processes}
\label{sec:br-proc}

When $\alpha=0$, $\mathbf{P}^K(b,d,0,z)$ is the law of a binary
branching process divided by $K$. Let us give some results on these
processes.
\begin{thm}
  \label{thm:invasion-extinction}
  Let $b,d>0$. As in Theorem~\ref{thm:exit-time}, define, for any
  $K\geq 1$ and any $z\in\mathbb{N}/K$,
  $\mathbf{P}^K_z=\mathbf{P}^K(b,d,0,z)$. Define also, for any
  $\rho\in\mathbb{R}$, on $\mathbb{D}(\mathbb{R}_+,\mathbb{R})$, the
  stopping time
  \begin{equation*}
    T_{\rho}=\inf\{t\geq 0:w_t=\rho\}.
  \end{equation*}
  Finally, let $(t_K)_{K\geq 1}$ be a sequence of positive numbers
  such that $\:\log K\ll t_K$.
  \begin{description}
  \item[\textmd{(a)}] If $b<d$ (sub-critical case), for any $\varepsilon>0$,
    \begin{gather}
      \label{eq:non-inv-subcrit}
      \lim_{K\rightarrow+\infty}\mathbf{P}^K_{1/K}(T_0\leq t_K\wedge
      T_{\lceil\varepsilon K\rceil/K})=1, \\
      \mbox{and}\quad \label{eq:ext-subcrit}
      \lim_{K\rightarrow+\infty}\mathbf{P}^K_{\lfloor\varepsilon K\rfloor/K}
      (T_0\leq t_K)=1.
    \end{gather}
    Moreover, for any $K\geq 1$, $k\geq 1$ and $n\geq 1$,
    \begin{equation}
      \label{eq:ext-before-expl}
      \mathbf{P}^K_{n/K}(T_{kn/K}\leq T_0)\leq\frac{1}{k}.
    \end{equation}
  \item[\textmd{(b)}] If $b>d$ (super-critical case), for any
    $\varepsilon>0$,
    \begin{gather}
      \label{eq:ext-supercrit}
      \lim_{K\rightarrow+\infty}\mathbf{P}^K_{1/K}
      (T_0\leq t_K\wedge
      T_{\lceil\varepsilon K\rceil/K})=\frac{d}{b} \\
      \mbox{and}\quad \label{eq:inv-supercrit}
      \lim_{K\rightarrow+\infty}\mathbf{P}^K_{1/K}
      (T_{\lceil\varepsilon K\rceil/K}\leq t_K)=1-\frac{d}{b}.
    \end{gather}
  \end{description}
\end{thm}

\paragraph{Proof}
Let us denote by $\mathbf{Q}_n$ the law of the binary branching
process with initial state $n\in\mathbb{N}$, with individual birth
rate $b$ and individual death rate $d$.
Then~(\ref{eq:non-inv-subcrit}), (\ref{eq:ext-subcrit}),
(\ref{eq:ext-before-expl}), (\ref{eq:ext-supercrit})
and~(\ref{eq:inv-supercrit}) rewrite respectively
\begin{gather}
  \lim_{K\rightarrow+\infty}\mathbf{Q}_1(T_0\leq t_K\wedge
  T_{\lceil\varepsilon K\rceil})=1, \label{eq:pf-BP-1} \\
  \lim_{K\rightarrow+\infty}\mathbf{Q}_{\lfloor\varepsilon K\rfloor}
  (T_0\leq t_K)=1, \label{eq:pf-BP-2} \\
  \mathbf{Q}_n(T_{kn}\leq T_0)\leq\frac{1}{k}, \label{eq:pf-BP-3} \\
  \lim_{K\rightarrow+\infty}\mathbf{Q}_1(T_0\leq t_K\wedge
  T_{\lceil\varepsilon K\rceil})=\frac{d}{b} \label{eq:pf-BP-4} \\
  \mbox{and}\quad \lim_{K\rightarrow+\infty}\mathbf{Q}_1
  (T_{\lceil\varepsilon K\rceil}\leq t_K)=1-\frac{d}{b}.\label{eq:pf-BP-5} 
\end{gather}

The limit~(\ref{eq:pf-BP-2}) follows easily from the distribution of
the extinction time for binary branching processes when $b\not= d$
(cf.\ Athreya and Ney, 1972, p.\ 109): for any $t\geq 0$ and
$n\in\mathbb{N}$,
\begin{equation}
  \label{eq:br-thm-non-crit}
  \mathbf{Q}_n(T_0\leq t)=\left(\frac{d\left(1-e^{-(b-d)t}\right)}
    {b-de^{-(b-d)t}}\right)^n.
\end{equation}

Since $t_K\rightarrow+\infty$, $\mathbf{Q}_1(T_0\leq
t_K\wedge T_{\lceil\varepsilon
  K\rceil})\rightarrow\mathbf{Q}_1(T_0<\infty)$, which
gives~(\ref{eq:pf-BP-1}) and~(\ref{eq:pf-BP-4}) (the probability of
extinction of a binary branching process can be recovered easily
from~(\ref{eq:br-thm-non-crit})\:).

The inequality~(\ref{eq:pf-BP-3}) follows from the fact that, if
$(Z_t,t\geq 0)$ is a process with law $\mathbf{Q}_n$,
$(Z_t\exp(-(b-d)t),t\geq 0)$ is a martingale (cf.\ Athreya and Ney,
1972, p.\ 111). Then, Doob's stopping theorem applied to the stopping
time $T_0\wedge T_{kn}$ yields,
\begin{equation*}
  \mathbf{E}_n(kne^{(d-b)T_{kn}}\mathbf{1}_{\{T_{kn}<T_0\}})=n,  
\end{equation*}
where $\mathbf{E}_n$ is the expectation with respect to
$\mathbf{Q}_n$. Therefore, when $b<d$, $kn\mathbf{Q}_n(T_{kn}<T_0)\leq
n$, and the proof of~(\ref{eq:pf-BP-3}) is completed.

The limit~(\ref{eq:pf-BP-5}) follows from the fact that, if
$(Z_t,t\geq 0)$ is a branching process with law $\mathbf{Q}_1$, the
martingale $(Z_t\exp(-(b-d)t),t\geq 0)$ converges a.s.~when
$t\rightarrow +\infty$ to a random variable $W$, where $W=0$ on the
event $\{T_0<\infty\}$ and $W>0$ on the event $\{T_0=\infty\}$ (cf.\ 
Athreya and Ney, 1972, p.\ 112). Hence, on the event $\{T_0=\infty\}$,
when $b>d$,
\begin{equation*}
  (\log Z_t)/t\rightarrow b-d>0.
\end{equation*}
Therefore, since $\log K\ll t_K$, for any $\varepsilon>0$,
$\mathbf{Q}_1(T_0=\infty,\ T_{\lceil\varepsilon K\rceil}\geq
t_K)\rightarrow 0$ when $K\rightarrow+\infty$.
Then,~(\ref{eq:pf-BP-5}) follows from the fact that
$\mathbf{Q}_1(T_0=\infty)=1-d/b$.\hfill$\Box$

\section{Proof of Theorem~\ref{thm:IPS-TSS}}
\label{sec:proof}

Let us assume, without loss of generality, that $\nu^K$ is constructed
by~(\ref{eq:def-XK}) on a sufficiently large probability space
$(\Omega,{\cal F},\mathbf{P})$.

We introduce the following sequences of stopping times: for all $n\geq
1$, let $\tau_n$ be the first mutation time after time $\tau_{n-1}$,
with $\tau_0=0$ (i.e.\ $\tau_n$ is the $n^{\mbox{\footnotesize{th}}}$
mutation time), and for any $n\geq 0$, let $\theta_n$ be the first
time after $\tau_n$ when the population gets monomorphic. Observe that
$\theta_0=0$ if the initial population is monomorphic. For any $n\geq
1$, define the random variable $U_n$ as the new trait value appearing
at the mutation time $\tau_n$, and, when $\theta_n<\infty$, define
$V_n$ by $\mbox{Supp}(\nu^K_{\theta_n})=\{V_n\}$. When
$\theta_n=+\infty$, define $V_n=+\infty$.

Our proof of Theorem~\ref{thm:IPS-TSS} is based on the following two
lemmas. The first lemma proves that there is no accumulation of
mutations on the time scale of Theorem~\ref{thm:IPS-TSS}, and studies
the asymptotic behavior of $\tau_1$ starting from a monomorphic
population, when $K\rightarrow +\infty$.
\begin{lembr}
  \label{lem:tau}
  \begin{description}
  \item[\textmd{(a)}] Assume that the initial condition of $\nu^K$
    satisfies
    $\sup_K\mathbf{E}(\langle\nu^K_0,\mathbf{1}\rangle)<+\infty$.
    Then, for any $\eta>0$, there exists $\varepsilon>0$ such that,
    for any $t>0$,
    \begin{equation}
      \label{eq:lem-tau-(a)-1}
      \limsup_{K\rightarrow+\infty}\mathbf{P}^K_{\nu^K_0}
      \left(\exists n\geq 0:\frac{t}{Ku_K}\leq\tau_n\leq
        \frac{t+\varepsilon}{Ku_K}\right)<\eta.
    \end{equation}
  \end{description}
  \noindent Let $x\in{\cal X}$ and let $(z_K)_{K\geq 1}$ be a sequence of
  integers such that $z_K/K\rightarrow z>0$.
  \begin{description}
  \item[\textmd{(b)}] For any $\varepsilon>0$,
    \begin{equation}
      \label{eq:lem-tau-(b)}
      \lim_{K\rightarrow+\infty}\mathbf{P}^K_{\frac{z_K}{K}\delta_x}
      \left(\tau_1>\log K,\ \sup_{t\in[\log K,\tau_1]}|\langle
      \nu^K_t,\mathbf{1}\rangle-\bar{n}_x|>\varepsilon\right)=0.
    \end{equation}
    Since $\log K\ll 1/Ku_K$, by~(a) with $t=0$,
    \begin{equation*}
      \lim_{K\rightarrow +\infty}\mathbf{P}^K_{\frac{z_K}{K}\delta_x}
      (\tau_1<\log K)=0.
    \end{equation*}
    In particular, under $\mathbf{P}^K_{\frac{z_K}{K}\delta_x}$,
    $\nu^K_{\log K}\overset{\cal P}{\rightarrow}\bar{n}_x\delta_x$ and
    $\nu^K_{\tau_1-}\overset{\cal P}{\rightarrow}\bar{n}_x\delta_x$.

    If, moreover, $z=\bar{n}_x$, then, for any $\varepsilon>0$,
    \begin{equation}
      \label{eq:lem-tau-(b)-bar-n}
      \lim_{K\rightarrow+\infty}\mathbf{P}^K_{\frac{z_K}{K}\delta_x}
      \left(\sup_{t\in[0,\tau_1]}|\langle
      \nu^K_t,\mathbf{1}\rangle-\bar{n}_x|>\varepsilon\right)=0.
    \end{equation}
  \item[\textmd{(c)}] For any $t>0$,
    \begin{equation*}
      \lim_{K\rightarrow+\infty}\mathbf{P}^K_{\frac{z_K}{K}\delta_x}
      \left(\tau_1>\frac{t}{Ku_K}\right)=\exp(-\beta(x)t),
    \end{equation*}
    where $\beta(\cdot)$ has been defined in~(\ref{eq:def-beta}).
  \end{description}
\end{lembr}
The second lemma studies the asymptotic behavior of $\theta_0$ and
$V_0$ starting from a dimorphic population, when $K\rightarrow
+\infty$.
\begin{lemma}
  \label{lem:theta}
  Fix $x,y\in{\cal X}$ satisfying~(\ref{eq:hyp-B1})
  or~(\ref{eq:hyp-B2}), and let $(z_K)_{K\geq 1}$ be a sequence of
  integers such that $z_K/K\rightarrow\bar{n}_x$. Then,
  \begin{gather}
    \label{eq:lem-theta-1}
    \lim_{K\rightarrow+\infty}\mathbf{P}^K_{\frac{z_K}{K}\delta_x
      +\frac{1}{K}\delta_y}(V_0=y)=\frac{[f(y,x)]_+}{b(y)}, \\
    \label{eq:lem-theta-2}
    \lim_{K\rightarrow+\infty}\mathbf{P}^K_{\frac{z_K}{K}\delta_x
      +\frac{1}{K}\delta_y}(V_0=x)=1-\frac{[f(y,x)]_+}{b(y)}, \\
    \label{eq:lem-theta-3}
    \forall\eta>0,\quad
    \lim_{K\rightarrow+\infty}\mathbf{P}^K_{\frac{z_K}{K}\delta_x
      +\frac{1}{K}\delta_y}
    \left(\theta_0>\frac{\eta}{Ku_K}\wedge\tau_1\right)=0 \\
    \label{eq:lem-theta-4}
    \mbox{and}\quad\forall\varepsilon>0,\quad
    \lim_{K\rightarrow+\infty}\mathbf{P}^K_{\frac{z_K}{K}\delta_x
      +\frac{1}{K}\delta_y}\left(|\langle\nu^K_{\theta_0},
      \mathbf{1}\rangle-\bar{n}_{V_0}|<\varepsilon\right)=1,
  \end{gather}
  where $f(y,x)$ has been defined in~(\ref{eq:def-fitn}).
\end{lemma}
Observe that~(\ref{eq:lem-theta-3}) implies in particular that
\begin{equation*}
  \lim_{K\rightarrow +\infty}
  \mathbf{P}^K_{\frac{z_K}{K}\delta_x+\frac{1}{K}\delta_y}
  (\theta_0<\tau_1)=1.
\end{equation*}

The proofs of these lemmas are postponed at the end of this section.

\paragraph{Proof of Theorem~\ref{thm:IPS-TSS}}
Observe that the generator $A$, defined in~(\ref{eq:generator-TSS}),
of the TSS process $(X_t,t\geq 0)$ of Theorem~\ref{thm:IPS-TSS} can be
written as
\begin{equation}
  \label{eq:generator-TSS-mu(x,dh)}
  A\varphi(x)=\int_{\mathbb{R}^l}(\varphi(x+h)-\varphi(x))
  \beta(x)\kappa(x,dh),
\end{equation}
where the probability measure $\kappa(x,dh)$ is defined by
\begin{multline}
  \label{eq:def-kappa}
  \kappa(x,dh)=\left(1-\int_{\mathbb{R}^l}
    \frac{[f(x+v,x)]_+}{b(x+v)}m(x,v)dv\right)
  \delta_0(dh) \\ +\frac{[f(x+h,x)]_+}{b(x+h)}m(x,h)dh.
\end{multline}
This means that the TSS model $X$ with initial state $x$ can be
constructed as follows: let $(Z(k),k=0,1,2,\ldots)$ be a Markov chain
in ${\cal X}$ with initial state $x$ and with transition kernel
$\kappa(x,dh)$, and let $(N(t),t\geq 0)$ be an independent standard
Poisson process. Then, the process $(X_t,t\geq
0)$ defined by
\begin{equation*}
  X_t:=Z\left(N\left(\int_0^t\beta(X_s)ds\right)\right)  
\end{equation*}
is a Markov process with infinitesimal
generator~(\ref{eq:generator-TSS-mu(x,dh)}) (cf.\ Ethier and Kurtz,
1986, Chap.~6). Let $\mathbf{P}_x$ denote its law, let $(T_n)_{n\geq
  1}$ denote the sequence of jump times of the Poisson process $N$ and
define $(S_n)_{n\geq 1}$ by $T_n=\int_0^{S_n}\beta(X_s)ds$. By~(A1)
and~(A3), $\beta(\cdot)>0$, and so $S_n$ is finite for any $n\geq 1$.
Observe that any jump of the process $X$ occurs at some time $S_n$,
but that all $S_n$ may not be effective jump times for $X$, because of
the Dirac mass at 0 appearing in~(\ref{eq:def-kappa}).

Fix $t>0$, $x\in{\cal X}$ and a measurable subset $\Gamma$ of ${\cal
  X}$. Under $\mathbf{P}_x$, $S_1$ and $X_{S_1}$ are independent,
$S_1$ is an exponential random variable with parameter $\beta(x)$, and
$X_{S_1}-x$ has law $\kappa(x,\cdot)$.  Therefore, for any $n\geq 1$,
the strong Markov property applied to $X$ at time $S_1$ yields
\begin{multline}
  \label{eq:induction-X}
  \mathbf{P}_x(S_n\leq t<S_{n+1},\ X_t\in\Gamma) \\
  =\int_0^t\beta(x)e^{-\beta(x)s}\int_{\mathbb{R}^l}
  \mathbf{P}_{x+h}(S_{n-1}\leq t-s<S_n,\ X_{t-s}\in\Gamma)
  \kappa(x,dh)ds.
\end{multline}
Moreover,
\begin{equation}
  \label{eq:init-induction-X}
  \mathbf{P}_x(0\leq t<S_1,\ X_t\in\Gamma)
  =\mathbf{1}_{\{x\in\Gamma\}}e^{-\beta(x)t}.
\end{equation}

The idea of our proof of Theorem~\ref{thm:IPS-TSS} is to show that the
same relations hold when we replace $S_n$ by $\tau_n$ and $X_t$ by the
support of $\nu^K_{t/Ku_K}$ (when it is a singleton) and when
$K\rightarrow +\infty$.

More precisely, fix $x\in{\cal X}$, $t>0$ and a measurable subset
$\Gamma$ of ${\cal X}$, and observe that
\begin{equation}
  \label{eq:pf-1}
  \left\{\exists y\in\Gamma:\mbox{Supp}(\nu^K_{t/Ku_K})=\{y\},\
  |\langle\nu^K_{t/Ku_K},\mathbf{1}\rangle-\bar{n}_y|<\varepsilon\right\}
  =\bigcup_{n\geq 0}A_n^K(t,\Gamma,\varepsilon),
\end{equation}
where
\begin{equation*}
  A_n^K(t,\Gamma,\varepsilon):=\left\{\theta_n\leq
    \frac{t}{Ku_K}<\tau_{n+1},\ V_n\in\Gamma,\ |\langle
    \nu^K_{t/Ku_K},\mathbf{1}\rangle-\bar{n}_{V_n}|<\varepsilon\right\}.
\end{equation*}
Let us define, for any $z\in\mathbb{N}$ and $n\geq 0$,
\begin{multline*}
  p^K_n(t,x,\Gamma,\varepsilon,z):=\mathbf{P}^K_{\frac{z}{K}\delta_x}
  \left(\theta_n\leq\frac{t}{Ku_K}<\tau_{n+1},\ V_n\in\Gamma,\right.
  \\ \sup_{s\in[\theta_n,\tau_{n+1}]}|\langle
  \nu^K_s,\mathbf{1}\rangle-\bar{n}_{V_n}|<\varepsilon\biggr)
\end{multline*}
and
\begin{align*}
  q^K_0(t,x,\Gamma,\varepsilon,z) &
  :=\mathbf{P}^K_{\frac{z}{K}\delta_x} \left(\frac{t}{Ku_K}<\tau_1,\ 
    V_0\in\Gamma,\ \sup_{s\in[\log K,\tau_1]}|\langle
    \nu^K_s,\mathbf{1}\rangle-\bar{n}_{V_0}|<\varepsilon\right) \\
  & =\mathbf{1}_{\{x\in\Gamma\}}\mathbf{P}^K_{\frac{z}{K}\delta_x}
  \left(\frac{t}{Ku_K}<\tau_1,\ \sup_{s\in[\log K,\tau_1]}|\langle
    \nu^K_s,\mathbf{1}\rangle-\bar{n}_x|<\varepsilon\right).
\end{align*}
Let us also extend these definitions to $\varepsilon=\infty$ by
suppressing the condition involving the supremum of
$|\langle\nu^K,\mathbf{1}\rangle-\bar{n}_{V_n}|$.

Then
\begin{lembr}
  \label{lem:induction}
  \begin{description}
  \item[\textmd{(a)}] For any $x\in{\cal X}$, $n\geq 1$, $t>0$,
    $\varepsilon\in]0,\infty]$ and for any sequence of integers
    $(z_K)$ such that $z_K/K\rightarrow z>0$, $p_n(t,x,\Gamma):=
    \lim_{K\rightarrow+\infty}p^K_n(t,x,\Gamma,\varepsilon,z_K)$
    exists, and is independent of $(z_K)$, $z$ and $\varepsilon$.
    
    Similarly, $p_0(t,x,\Gamma):=
    \lim_{K\rightarrow+\infty}q^K_0(t,x,\Gamma,\varepsilon,z_K)$
    exists, and is independent of $(z_K)$, $z$ and $\varepsilon$,
    and, if $z=\bar{n}_x$,
    $\lim_{K\rightarrow+\infty}p^K_0(t,x,\Gamma,\varepsilon,z_K)$
    exists and is also equal to $p_0(t,x,\Gamma)$.
    
    Finally, if we assume that $(z_K)$ is a sequence of
    $\mathbb{N}$-valued random variables such that $z_K/K$ converge in
    probability to a deterministic $z>0$, then the limits above hold
    in probability (with the same restriction that $z$ has to be equal
    to $\bar{n}_x$ for $p^K_0$).
  \item[\textmd{(b)}] The functions $p_n(t,x,\Gamma)$ are continuous
    with respect to $t$ and measurable with respect to $x$, and
    satisfy
    \begin{gather}
      p_0(t,x,\Gamma)=\mathbf{1}_{\{x\in\Gamma\}}
      e^{-\beta(x)t} \quad\mbox{and}\quad\forall n\geq 0, \notag \\
      p_{n+1}(t,x,\Gamma)=\int_0^t\beta(x)e^{-\beta(x)s}
      \int_{\mathbb{R}^l}p_n(t-s,x+h,\Gamma)\kappa(x,dh)ds.
      \label{eq:induction-pn}
    \end{gather}
  \end{description}
\end{lembr}
Let us postpone the proof of this lemma after the proof of
Theorem~\ref{thm:IPS-TSS}.

Observe that, because of~(\ref{eq:induction-X})
and~(\ref{eq:init-induction-X}), Lemma~\ref{lem:induction}~(b) implies
that $\mathbf{P}_x(S_n\leq t<S_{n+1},\ X_t\in\Gamma)=p_n(t,x,\Gamma)$.

Now, let $\mathbf{\tilde{P}}^K_{\nu}$ denote the law of the process
$\nu^K$ with random initial state $\nu$. Since $\nu^K$ is Markov,
$\mathbf{\tilde{P}}^K_{\gamma_K/K\delta_x}
=\mathbf{E}[\mathbf{P}^K_{\gamma_K(\omega)/K\delta_x}]$.
By~(\ref{eq:pf-1}),
\begin{multline*}
  \mathbf{\tilde{P}}^K_{\frac{\gamma_K}{K}\delta_x}\left(\exists
    y\in\Gamma:\mbox{Supp} (\nu^K_{t/Ku_K})=\{y\},\ \right. \\
  \left. |\langle
    \nu^K_{t/Ku_K},\mathbf{1}\rangle-\bar{n}_y|<\varepsilon\right)
  =\sum_{n\geq 0}\mathbf{\tilde{P}}^K_{\frac{\gamma_K}{K}\delta_x}
  (A^K_n(t,\Gamma,\varepsilon)),
\end{multline*}
where $(\gamma_K)$ is the sequence of $\mathbb{N}$-valued random
variables of Theorem~\ref{thm:IPS-TSS}.

For any $K\geq 1$ and $n\geq 1$,
\begin{gather*}
  p_n^K(t,x,\Gamma,\varepsilon,\gamma_K)
  \leq\mathbf{P}^K_{\frac{\gamma_K}{K}\delta_x}
  (A^K_n(t,\Gamma,\varepsilon))
  \leq p_n^K(t,x,\Gamma,\infty,\gamma_K), \\
  \mbox{and}\quad q_0^K(t,x,\Gamma,\varepsilon,\gamma_K)
  \leq\mathbf{P}^K_{\frac{\gamma_K}{K}\delta_x}
  (A^K_n(t,\Gamma,\varepsilon))
  \leq p_n^K(t,x,\Gamma,\infty,\gamma_K),
\end{gather*}
so, by Lemma~\ref{lem:induction}~(a), for any $n\geq 0$,
$\mathbf{P}^K_{(\gamma_K/K)\delta_x}
(A^K_n(t,\Gamma,\varepsilon))\overset{\cal
  P}{\rightarrow}p_n(t,x,\Gamma)$, and therefore,
\begin{equation}
  \label{eq:lim-Ptilde-pn}
  \lim_{K\rightarrow+\infty}\mathbf{\tilde{P}}^K_{(\gamma_K/K)\delta_x}
  (A^K_n(t,\Gamma,\varepsilon))=p_n(t,x,\Gamma).
\end{equation}

Now, by~(\ref{eq:pf-1}), for any $K\geq 1$,
\begin{equation}
  \label{eq:bound-Ptilde}
  \sum_{n=0}^{+\infty}\left[\mathbf{\tilde{P}}^K_{\frac{\gamma_K}{K}\delta_x}
  (A^K_n(t,\Gamma,\varepsilon))
  +\mathbf{\tilde{P}}^K_{\frac{\gamma_K}{K}\delta_x}
  (A^K_n(t,\Gamma^c,\varepsilon))\right]\leq 1,
\end{equation}
where $\Gamma^c$ denotes the complement of $\Gamma$. Moreover,
$\sum_{n=0}^{+\infty}[p_n(t,x,\Gamma)+p_n(t,x,\Gamma^c)]=1$.
Therefore, for any $\eta>0$, there exists $n_0$ such that
\begin{equation*}
  \sum_{n=0}^{n_0}[p_n(t,x,\Gamma)+p_n(t,x,\Gamma^c)]\geq 1-\eta.
\end{equation*}
Then, one can easily deduce from~(\ref{eq:lim-Ptilde-pn})
and~(\ref{eq:bound-Ptilde}) that
\begin{equation*}
  \limsup_{K\rightarrow+\infty}\sum_{n\geq n_0}
  \left[\mathbf{\tilde{P}}^K_{\frac{\gamma_K}{K}\delta_x}
  (A^K_n(t,\Gamma,\varepsilon))
  +\mathbf{\tilde{P}}^K_{\frac{\gamma_K}{K}\delta_x}
  (A^K_n(t,\Gamma^c,\varepsilon))\right]\leq \eta,
\end{equation*}
from which follows, by~(\ref{eq:pf-1}), that
\begin{multline*}
  \lim_{K\rightarrow+\infty}
  \mathbf{\tilde{P}}^K_{\frac{\gamma_K}{K}\delta_x}\left(
    \exists y\in\Gamma:\mbox{Supp}
    (\nu^K_{t/Ku_K})=\{y\},\ y\in\Gamma,\ |\langle
    \nu^K_{t/Ku_K},\mathbf{1}\rangle-\bar{n}_y|<\varepsilon\right) \\
  =\sum_{n\geq 0}p_n(t,x,\Gamma)=\mathbf{P}_x(X_t\in\Gamma),
\end{multline*}
which is~(\ref{eq:def-cvgce}) in the case of a single time $t$.

In order to complete the proof of Theorem~\ref{thm:IPS-TSS}, we have
to generalize this limit to any sequence of times $0<t_1<\ldots<t_n$.

We will specify the method only in the case of two times $0<t_1<t_2$.
It can be easily generalized to a sequence of $n$ times. We introduce
for any integers $0\leq n_1\leq n_2$ the probabilities
\begin{multline*}
    p^K_{n_1,n_2}(t_1,t_2,x,\Gamma_1,\Gamma_2,\varepsilon,z)
    \\ :=\mathbf{P}^K_{\frac{z}{K}\delta_x}
    \left(\theta_{n_1}\leq\frac{t_1}{Ku_K}<\tau_{n_1+1},\
      V_{n_1}\in\Gamma_1,\!\sup_{s\in[\theta_{n_1},\tau_{n_1+1}]}|\langle
      \nu^K_s,\mathbf{1}\rangle-\bar{n}_{V_{n_1}}|<\varepsilon,\ 
    \right. \\
    \left. \theta_{n_2}\leq\frac{t_2}{Ku_K}<\tau_{n_2+1},\ 
      V_{n_2}\in\Gamma_2\mbox{\ and\
    }\sup_{s\in[\theta_{n_2},\tau_{n_2+1}]}|\langle
      \nu^K_s,\mathbf{1}\rangle-\bar{n}_{V_{n_2}}|<\varepsilon\right),
\end{multline*}
and
\begin{multline*}
    q^K_{0,n_2}(t_1,t_2,x,\Gamma_1,\Gamma_2,\varepsilon,z)
    \\ :=\mathbf{1}_{\{x\in\Gamma_1\}}\mathbf{P}^K_{\frac{z}{K}\delta_x}
    \left(\frac{t_1}{Ku_K}<\tau_1,\
      \sup_{s\in[\log K,\tau_1]}|\langle
      \nu^K_s,\mathbf{1}\rangle-\bar{n}_x|<\varepsilon,\ 
    \right. \\
    \left. \theta_{n_2}\leq\frac{t_2}{Ku_K}<\tau_{n_2+1},\ 
      V_{n_2}\in\Gamma_2\mbox{\ and\
    }\sup_{s\in[\theta_{n_2},\tau_{n_2+1}]}|\langle
      \nu^K_s,\mathbf{1}\rangle-\bar{n}_{V_{n_2}}|<\varepsilon\right).
\end{multline*}

Then, we can use a calculation very similar to the proof of
Lemma~\ref{lem:induction} to prove that, as $K\rightarrow +\infty$,
$p^K_{n_1,n_2}(t_1,t_2,x,\Gamma_1,\Gamma_2,\varepsilon,z_K)$ converges
to a limit $p_{n_1,n_2}(t_1,t_2,x,\Gamma_1,\Gamma_2)$ independent of
$\varepsilon\in]0,\infty]$, $z_K$ and the limit $z>0$ of $z_K/K$ (with
the restriction that $z$ has to be equal to $\bar{n}_x$ if $n_1=0$),
and that $\lim q^K_{0,n_2}(t_1,t_2,x,\Gamma_1,\Gamma_2,\varepsilon,z)
=p_{0,n_2}(t_1,t_2,x,\Gamma_1,\Gamma_2)$,
where
\begin{equation*}
  \left\{\begin{array}{l}
      \displaystyle{p_{0,n_2}(t_1,t_2,x,\Gamma_1,\Gamma_2)
        =\mathbf{1}_{\{x\in\Gamma_1\}}
        e^{-\beta(x)t_1}p_{n_2}(t_2-t_1,x,\Gamma_2)}; \\
      \displaystyle{p_{n_1+1,n_2+1}(t_1,t_2,x,\Gamma_1,\Gamma_2)
        \phantom{\int}} \\ \qquad
        \displaystyle{=\int_0^{t_1}\beta(x)e^{-\beta(x)s}
        \int_{\mathbb{R}^l}p_{n_1,n_2}(t_1-s,t_2-s,x+h,\Gamma_1,\Gamma_2)
        \kappa(x,dh)ds.}
    \end{array}\right.
\end{equation*}
As above, we obtain equation~(\ref{eq:def-cvgce}) for $n=2$ by
observing that the same relation holds for the TSS process $X$.

This completes the proof of Theorem~\ref{thm:IPS-TSS}.\hfill$\Box$

\paragraph{Proof of Lemma~\ref{lem:induction}}
First, let us prove that the convergence of
$p_n^K(t,x,\Gamma,\varepsilon,z_K)$ when $z_K\in\mathbb{N}$ in
Lemma~\ref{lem:induction}~(a) implies the convergence in probability
of these quantities when $z_K$ are random variables: if $(z_K)$ is a
sequence of random variables such that $z_K/K\overset{\cal
  P}{\rightarrow} z$, by Skorohod's Theorem, we can construct on an
auxiliary probability space $\hat{\Omega}$ a sequence of random
variables $(\hat{z}_K)$ such that ${\cal L}(\hat{z}_K)={\cal L}(z_K)$
and $\hat{z}_K(\hat{\omega})/K\rightarrow z$ for any
$\hat{\omega}\in\hat{\Omega}$. Then, $\lim
p^K_n(t,x,\Gamma,\varepsilon,\hat{z}_K(\hat{\omega}))=p_n(t,x,\Gamma)$
for any $\hat{\omega}\in\hat{\Omega}$, which implies that
$p^K_n(t,x,\Gamma,\varepsilon,z_K)\overset{\cal
  P}{\rightarrow}p_n(t,x,\Gamma)$. The same method applies to
$q^K_0(t,x,\Gamma,\varepsilon,z_k)$.

We will prove Lemma~\ref{lem:induction}~(a) and~(b) by induction over
$n\geq 0$.

First, when $t>0$, it follows from the fact that $t/Ku_K>\log K$ for
sufficiently large $K$, and from Lemma~\ref{lem:tau}~(b) and~(c), that
\begin{equation*}
  \lim_{K\rightarrow+\infty}q_0^K(t,x,\Gamma,\varepsilon,z_K)
  =\mathbf{1}_{\{x\in\Gamma\}}e^{-\beta(x)t},
\end{equation*}
and that, if $z=\bar{n}_x$,
\begin{equation*}
  \lim_{K\rightarrow+\infty}p_0^K(t,x,\Gamma,\varepsilon,z_K)
  =\mathbf{1}_{\{x\in\Gamma\}}e^{-\beta(x)t}.
\end{equation*}

Then, fix $n\geq 0$ and assume that Lemma~\ref{lem:induction}~(a)
holds for $n$. We intend to prove the convergence of
$p^K_{n+1}(t,x,\Gamma,\varepsilon,z_K)$ to $p_{n+1}(t,x,\Gamma)$
satisfying~(\ref{eq:induction-pn}) by applying the strong Markov
property at time $\tau_1$, in a similar way as when we
obtained~(\ref{eq:induction-X}). However, the convergence of
$p^K_n(t,x,\Gamma,\varepsilon,z_K)$ to $p_n(t,x,\Gamma)$ only holds
for \emph{non-random} $t$. Therefore, we will divide the time interval
$[0,t]$ in a finite number of small intervals and use the Markov
property at time $\tau_1$ when $\tau_1$ is in each of these intervals.
Moreover, we will also use the Markov property at time $\theta_1$ and
we will use the fact that $U_1$ is independent of $\tau_1$ and
$\nu^K_{\tau_1-}$ and that $U_1-x$ is a random variable with law
$m(x,h)dh$.

Following this program, we can bound
$p^K_{n+1}(t,x,\Gamma,\varepsilon,z_K)$ from above as follows: fix
$\eta>0$; using Lemma~\ref{lem:tau}~(a) in the first inequality, for
sufficiently large $k\geq 0$ and $K\geq 1$,
\begin{align}
  & p^K_{n+1}(t,x,\Gamma,\varepsilon,z_K)\leq
  \mathbf{P}^K_{\frac{z_K}{K}\delta_x}
  \left(\theta_{n+1}\leq\frac{t}{Ku_K},\ 
    \tau_{n+2}>\frac{t+2/2^k}{Ku_K},\ 
    V_{n+1}\in\Gamma \right)+\eta \notag \\
  & \leq\sum_{i=0}^{\lceil t2^k\rceil-1}
  \mathbf{P}^K_{\frac{z_K}{K}\delta_x}\left(\frac{i}{2^kKu_K}
    \leq\tau_1\leq\frac{i+1}{2^kKu_K},\ 
    \theta_{n+1}\leq\frac{t}{Ku_K}, \right. \notag \\
  & \left.\qquad \tau_{n+2}>\frac{t+2/2^k}{Ku_K}\mbox{\ and\ }
    V_{n+1}\in\Gamma\right)+\eta \notag \\
  & \leq\sum_{i=0}^{\lceil t2^k\rceil-1}
  \mathbf{E}^K_{\frac{z_K}{K}\delta_x}\biggl[
    \mathbf{1}_{\left\{\frac{i}{2^kKu_K}\leq\tau_1
        \leq\frac{i+1}{2^kKu_K}\right\}}
    \mathbf{P}^K_{\nu^K_{\tau_1-}+\frac{1}{K}\delta_{U_1}}
    \left(\theta_n\leq\frac{t-i/2^k}{Ku_K}, \right. \notag \\
  & \left.\left.\qquad \tau_{n+1}>\frac{t-(i-1)/2^k}{Ku_K}\mbox{\ 
        and\ }V_n\in\Gamma\right)\right]+\eta \notag \\
  & \leq\sum_{i=0}^{\lceil t2^k\rceil-1}
  \mathbf{E}^K_{\frac{z_K}{K}\delta_x}\biggl[
    \mathbf{1}_{\left\{\frac{i}{2^kKu_K}\leq\tau_1
        \leq\frac{i+1}{2^kKu_K}\right\}}
    \int_{\mathbb{R}^l}\mathbf{E}^K_{\nu^K_{\tau_1-}+\frac{1}{K}\delta_{x+h}}
    \biggl(\mathbf{1}_{\left\{\theta_0\geq\frac{1}{2^kKu_K}
          \wedge\tau_1\right\}} \notag \\
  & \quad +\mathbf{1}_{\left\{\theta_0<\frac{1}{2^kKu_K}\wedge\tau_1\right\}}
      \mathbf{P}^K_{\nu^K_{\theta_0}}\left(\theta_n\leq\frac{t-i/2^k}{Ku_K}
        <\tau_{n+1},\ V_n\in\Gamma\right)\biggr)m(x,h)dh\biggr]+\eta.
      \notag \\
  & \leq\sum_{i=0}^{\lceil t2^k\rceil-1}
  \mathbf{E}^K_{\frac{z_K}{K}\delta_x}\biggl[
    \mathbf{1}_{\left\{\frac{i}{2^kKu_K}\leq\tau_1
        \leq\frac{i+1}{2^kKu_K}\right\}}
    \int_{\mathbb{R}^l}\mathbf{E}^K_{\nu^K_{\tau_1-}+\frac{1}{K}\delta_{x+h}}
    \biggl(\mathbf{1}_{\left\{\theta_0\geq\frac{1}{2^kKu_K}
          \wedge\tau_1\right\}} \notag \\
  & \quad +\mathbf{1}_{\left\{\theta_0<\frac{1}{2^kKu_K}\wedge\tau_1\right\}}
      p_n^K(t-i/2^k,V_0,\Gamma,\infty,K\langle
      \nu^K_{\theta_0},\mathbf{1}\rangle)\biggr)m(x,h)dh\biggr]+\eta.
  \label{eq:UB-pn-1}
\end{align}


Now, since
$\nu^K_{\tau_1-}=\langle\nu^K_{\tau_1-},\mathbf{1}\rangle\delta_x$,
under $\mathbf{P}^K_{\nu_{\tau_1-}^K+\frac{1}{K}\delta_{x+h}}$, on the
event $\{\theta_0<\tau_1\}$,
\begin{multline}
  \label{eq:UB-pn-2}
  p_n^K(t-i/2^k,V_0,\Gamma,\infty,K\langle
  \nu^K_{\theta_0},\mathbf{1}\rangle)=\mathbf{1}_{\{V_0=x\}}
  p_n^K(t-i/2^k,x,\Gamma,\infty,K\langle
  \nu^K_{\theta_0},\mathbf{1}\rangle) \\
  +\mathbf{1}_{\{V_0=x+h\}}
  p_n^K(t-i/2^k,x+h,\Gamma,\infty,K\langle
  \nu^K_{\theta_0},\mathbf{1}\rangle).
\end{multline}

By Lemma~\ref{lem:tau}~(b), $\nu^K_{\tau_1-}\overset{\cal
  P}{\rightarrow}\bar{n}_x\delta_x$ under
$\mathbf{P}_{\frac{z_K}{K}\delta_x}$, so we can use Skorohod's Theorem
to construct random variables $\hat{N}_K$ on an auxiliary probability
space $\hat{\Omega}$ with the same law that $\langle
\nu^K_{\tau_1-},\mathbf{1}\rangle$ and converging to $\bar{n}_x$ for
any $\hat{\omega}\in\hat{\Omega}$.

Fix $\hat{\omega}\in\hat{\Omega}$. Under
$\mathbf{P}^K_{\hat{N}_K(\hat{\omega})\delta_x+\frac{1}{K}\delta_{x+h}}$,
define
\begin{equation*}
  Z_1^K=\langle\nu^K_{\theta_0},\mathbf{1}\rangle\mathbf{1}_{\{V_0=x,\
    \theta_0<\tau_1\}}+\frac{\lceil K\bar{n}_x\rceil}{K}
  \mathbf{1}_{\{V_0\not=x\}\cup\{\theta_0\geq\tau_1\}}.
\end{equation*}
It follows from Lemma~\ref{lem:theta}~(\ref{eq:lem-theta-3})
and~(\ref{eq:lem-theta-4}), and from assumption~(B) that, for Lebesgue
almost every $h$, $Z^K_1\overset{\cal P}{\rightarrow}\bar{n}_x$, so,
by the induction assumption, under
$\mathbf{P}^K_{\hat{N}_K(\hat{\omega})\delta_x+\frac{1}{K}\delta_{x+h}}$,
\begin{equation*}
  p^K_n(t-i/2^k,x,\Gamma,\infty,KZ^K_1)
  \overset{\cal P}{\rightarrow}p_n(t-i/2^k,x,\Gamma).
\end{equation*}

Now, given two sequences of uniformly bounded random variables
$(X_K)_{K\geq 1}$ and $(Y_K)_{K\geq 1}$ such that $X_K$ and $Y_K$ are
defined on the same probability space for any $K\geq 1$, and such
that, when $K\rightarrow +\infty$, $X_K$ converges in probability to a
constant $C$ and $\lim_K\mathbf{E}(Y_K)$ exists, it is standard to
prove that
\begin{equation}
  \label{eq:cv-proba-bounded}
  \lim_{K\rightarrow +\infty}\mathbf{E}(X_KY_K)=C\lim_{K\rightarrow
  +\infty}\mathbf{E}(Y_K).
\end{equation}

Applying this with $X_K=p^K_n(t-i/2^k,x,\Gamma,\infty,KZ^K_1)$ and
$Y_K=\mathbf{1}_{\{V_0=x,\ \theta_0<\tau_1\}}$, by
Lemma~\ref{lem:theta}~(\ref{eq:lem-theta-2})
and~(\ref{eq:lem-theta-3}) and assumption~(B), for Lebesgue almost any
$h$, and for any $\hat{\omega}\in\hat{\Omega}$,
\begin{multline*}
  \lim_{K\rightarrow+\infty}
  \mathbf{E}^K_{\hat{N}_K(\hat{\omega})\delta_x+\frac{1}{K}\delta_{x+h}}
  \left(\mathbf{1}_{\{V_0=x,\ \theta_0<\tau_1\}}
    p_n^K(t-i/2^k,x,\Gamma,\infty,K\langle
    \nu^K_{\theta_0},\mathbf{1}\rangle)\right) \\
  =\left(1-\frac{[f(x+h,x)]_+}{b(x+h)}\right)p_n(t-i/2^k,x,\Gamma).
\end{multline*}
Finally, we obtain that, for Lebesgue almost any $h$, under
$\mathbf{P}^K_{\frac{z_K}{K}\delta_x}$,
\begin{multline}
  \label{eq:UB-pn-3}
  \mathbf{E}^K_{\nu^K_{\tau_1-}+\frac{1}{K}\delta_{x+h}}
  \left(\mathbf{1}_{\{V_0=x,\ \theta_0<\tau_1\}}
    p_n^K(t-i/2^k,x,\Gamma,\infty,K\langle
    \nu^K_{\theta_0},\mathbf{1}\rangle)\right) \\
  \overset{\cal P}{\rightarrow}
  \left(1-\frac{[f(x+h,x)]_+}{b(x+h)}\right)p_n(t-i/2^k,x,\Gamma).
\end{multline}

Similarly, we can use Lemma~\ref{lem:theta}~(\ref{eq:lem-theta-1}) and
the random variable
\begin{equation*}
  Z_2^K=\langle\nu^K_{\theta_0},\mathbf{1}\rangle\mathbf{1}_{\{V_0=x+h,\
    \theta_0<\tau_1\}}+\bar{n}_{x+h}
  \mathbf{1}_{\{V_0\not=x+h\}\cup\{\theta_0\geq\tau_1\}}
\end{equation*}
to prove that, for Lebesgue almost any $h$, under
$\mathbf{P}^K_{\frac{z_K}{K}\delta_x}$,
\begin{multline}
  \label{eq:UB-pn-4}
  \mathbf{E}^K_{\nu^K_{\tau_1-}+\frac{1}{K}\delta_{x+h}}\left(
  \mathbf{1}_{\{V_0=x+h,\ \theta_0<\tau_1\}}
  p_n^K(t-i/2^k,x+h,\Gamma,\infty,K\langle
  \nu^K_{\theta_0},\mathbf{1}\rangle)\right) \\
  \overset{\cal P}{\rightarrow}
  \frac{[f(x+h,x)]_+}{b(x+h)}p_n(t-i/2^k,x+h,\Gamma).
\end{multline}
Moreover, by Lemma~\ref{lem:theta}~(\ref{eq:lem-theta-3}), for
Lebesgue almost any $h$, under $\mathbf{P}^K_{(z_K/K)\delta_x}$,
\begin{equation}
  \label{eq:UB-pn-5}
  \mathbf{P}^K_{\nu^K_{\tau_1-}+\frac{1}{K}\delta_{x+h}}\left(
    \theta_0\geq\frac{1}{2^kKu_K}\wedge\tau_1\right)
  \overset{\cal P}{\rightarrow}0.  
\end{equation}

Collecting these results together,
applying~(\ref{eq:cv-proba-bounded}) again, it follows from
Lem\-ma~\ref{lem:tau}~(c) and~(\ref{eq:UB-pn-2})
that, for Lebesgue almost any $h$,
\begin{align*}
  & \lim_{K\rightarrow +\infty}
  \mathbf{E}^K_{\frac{z_K}{K}\delta_x}
  \biggl[\mathbf{1}_{\left\{\frac{i}{2^kKu_K}\leq\tau_1
      \leq\frac{i+1}{2^kKu_K}\right\}}
  \mathbf{E}^K_{\nu^K_{\tau_1-}+\frac{1}{K}\delta_{x+h}}
  \biggl(\mathbf{1}_{\left\{\theta_0\geq\frac{1}{2^kKu_K}
      \wedge\tau_1\right\}} \\
  & \qquad +\mathbf{1}_{\left\{\theta_0<\frac{1}
      {2^kKu_K}\wedge\tau_1\right\}}
  p_n^K(t-i/2^k,V_0,\Gamma,\infty,K\langle
  \nu^K_{\theta_0},\mathbf{1}\rangle)
  \biggr)\biggr] \\
  & =\left(e^{-\beta(x)\frac{i}{2^k}}-e^{-\beta(x)\frac{i+1}{2^k}}\right)
  \left[\frac{[f(x+h,x)]_+}{b(x+h)}p_n(t-i/2^k,x+h,\Gamma) \right. \\
  & \qquad\qquad \left.
    +\left(1-\frac{[f(x+h,x)]_+}{b(x+h)}\right)p_n(t-i/2^k,x,\Gamma)\right].
\end{align*}

Finally, taking the integral of both sides with respect to $m(x,h)dh$,
the dominated convergence theorem and~(\ref{eq:UB-pn-1}) yield
\begin{multline*}
  \limsup_{K\rightarrow+\infty}\ 
  p_{n+1}^K(x,t,\Gamma,\varepsilon,z_K) \\ \leq\sum_{i=0}^{\lceil
    t2^k\rceil-1}
  \left(e^{-\beta(x)\frac{i}{2^k}}-e^{-\beta(x)\frac{i+1}{2^k}}\right)
  \int_{\mathbb{R}^l}p_n(t-i/2^k,x+h,\Gamma)\kappa(x,dh)+\eta.
\end{multline*}
Taking the limit $k\rightarrow+\infty$ first and then $\eta\rightarrow
0$, it follows from the fact that
\begin{equation*}
  e^{-\beta(x)i/2^k}-e^{-\beta(x)(i+1)/2^k}=e^{-\beta(x)
  i/2^k}(\beta(x)/2^k+O(1/2^{2k}))
\end{equation*}
and from the convergence of Riemann sums that
\begin{equation*}
  \limsup_{K\rightarrow+\infty}\ p_{n+1}^K(x,t,\Gamma,\varepsilon,z_K) 
  \leq\int_0^t\beta(x)e^{-\beta(x)s}
  \int_{\mathbb{R}^l}p_n(t-s,x+h,\Gamma)\kappa(x,dh)ds. 
\end{equation*}

Using the same method as for~(\ref{eq:UB-pn-1}), we can give a lower
bound for $p^K_n$ as follows: for any $\eta>0$, for sufficiently large
$k\geq 0$ and $K\geq 1$,
\begin{align*}
  & p^K_{n+1}(t,x,\Gamma,\varepsilon,z_K)\geq
  \mathbf{P}^K_{\frac{z_K}{K}\delta_x}
  \left(\theta_{n+1}\leq\frac{t}{Ku_K},\
    \tau_{n+2}>\frac{t-2/2^k}{Ku_K},\
    V_{n+1}\in\Gamma \right. \\
  & \qquad \mbox{and\ }\sup_{s\in[\theta_{n+1},\tau_{n+2}]}|\langle
  \nu^K_s,\mathbf{1}\rangle-\bar{n}_{V_{n+1}}|<\varepsilon
  \biggr)-\eta \\
  & \geq\sum_{i=0}^{\lfloor t2^k\rfloor-3}
  \mathbf{E}^K_{\frac{z_K}{K}\delta_x}\biggl[
  \mathbf{1}_{\left\{\frac{i}{2^kKu_K}\leq\tau_1
      \leq\frac{i+1}{2^kKu_K}\right\}}
  \mathbf{P}^K_{\nu^K_{\tau_1-}+\frac{1}{K}\delta_{U_1}}
  \left(\theta_n\leq\frac{t-(i+1)/2^k}{Ku_K}, \right. \\
  & \qquad \tau_{n+1}>\frac{t-(i+2)/2^k}{Ku_K},\ 
  V_n\in\Gamma\mbox{\ and\ }
  \sup_{s\in[\theta_n,\tau_{n+1}]}|\langle
  \nu^K_s,\mathbf{1}\rangle-\bar{n}_{V_n}|<\varepsilon
  \biggr)\biggr]-\eta \\
  & \geq\sum_{i=0}^{\lfloor t2^k\rfloor-3}
  \mathbf{E}^K_{\frac{z_K}{K}\delta_x}\biggl[
    \mathbf{1}_{\left\{\frac{i}{2^kKu_K}\leq\tau_1
        \leq\frac{i+1}{2^kKu_K}\right\}}
    \int_{\mathbb{R}^l}\mathbf{E}^K_{\nu^K_{\tau_1-}+\frac{1}{K}\delta_{x+h}}
    \biggl(\mathbf{1}_{\left\{\theta_0<\frac{1}{2^kKu_K}
          \wedge\tau_1\right\}} \\
  & \qquad p_n^K(t-(i+2)/2^k,V_0,\Gamma,\varepsilon,K\langle
    \nu^K_{\theta_0},\mathbf{1}\rangle)\biggr)m(x,h)dh\biggr]-\eta.
\end{align*}

Then, as above, letting $K\rightarrow+\infty$, then
$k\rightarrow+\infty$ and finally $\eta\rightarrow 0$, we obtain
\begin{equation*}
  \liminf_{K\rightarrow+\infty}p_{n+1}^K(x,t,\Gamma,\varepsilon,z_K) 
  \geq\int_0^t\beta(x)e^{-\beta(x)s}
        \int_{\mathbb{R}^l}p_n(t-s,x+h,\Gamma)\kappa(x,dh)ds,
\end{equation*}
which completes the proof of Lemma~\ref{lem:induction} by
induction.\hfill$\Box$

\paragraph{Proof of Lemma~\ref{lem:tau}~(a)}
Fix $\eta>0$. By Theorem~\ref{thm:compar}~(a) and~(c), for any $K\geq
1$,
\begin{gather*}
  \langle\nu^K,\mathbf{1}\rangle\preceq Z^K, \\
  \mbox{where}\quad {\cal L}(Z^K)=\mathbf{P}^K(2\bar{b},0,\underline{\alpha},
  \langle\nu^K_0,\mathbf{1}\rangle+1).
\end{gather*}
Since $\sup_K\mathbf{E}(\langle\nu^K_0,\mathbf{1}\rangle)<+\infty$, we
can choose $M<+\infty$ such that
\begin{equation*}
  \sup_{K\geq 1}\mathbf{P}(\langle\nu^K_0,\mathbf{1}\rangle+1>M)<\eta/3.  
\end{equation*}

Then, apply Theorem~\ref{thm:exit-time}~(c) to
$\mathbf{P}^K(2\bar{b},0,\underline{\alpha},
\langle\nu^K_0,\mathbf{1}\rangle+1)$ with $C=[1,M]$, $\eta_2=M$ and
$\eta_1$ such that $0<2\bar{b}/\underline{\alpha}-\eta_1<1/2$: there
exists $V>0$ such that
\begin{gather}
  \label{eq:pf-lem-tau-(a)}
  \limsup_{K\rightarrow+\infty}\mathbf{P}(T^K<e^{KV})<\eta/3, \\
  \mbox{where}\quad T^K=\inf\{t\geq
  0,Z^K_t\not\in[1/2,M+2\bar{b}/\underline{\alpha}]\}. \notag
\end{gather}

Fix $t,\varepsilon>0$. Since, for $s\leq T^K$,
$\langle\nu^K_s,\mathbf{1}\rangle\leq M+2\bar{b}/\underline{\alpha}$,
if we apply Theorem~\ref{thm:compar}~(b) to the process
$(\nu^K_{s+(t/Ku_K)}-\nu^K_{t/Ku_K},\ s\geq 0)$, we obtain, for $s\leq
T^K-t/Ku_K$,
\begin{equation*}
  A^K_{s+(t/Ku_K)}-A^K_{t/Ku_K}\preceq B^K_s,
\end{equation*}
where $A^K_s$ is the number of mutations occuring between 0 and $s$,
and where $B^K$ is a Poisson process with parameter
$Ku_K\bar{b}(M+2\bar{b}/\underline{\alpha})$. Therefore,
combining~(\ref{eq:pf-lem-tau-(a)}) with the fact that $1/Ku_K\ll
e^{KV}$, we obtain that, for sufficiently large $K$
\begin{align*}
  \mathbf{P}(A^K_{(t+\varepsilon)/Ku_K}
  -A^K_{t/Ku_K}\geq 1) & \leq\mathbf{P}(B^K_{\varepsilon/Ku_K}\geq
  1)+2\eta/3 \\ & =1-\exp(-\bar{b}(M+2\bar{b}/\underline{\alpha})
  \varepsilon)+2\eta/3,
\end{align*}
which can be made smaller than $\eta$ if $\varepsilon$ is sufficiently
small. This ends the proof of~(\ref{eq:lem-tau-(a)-1}).\hfill$\Box$

\paragraph{Proof of Lemma~\ref{lem:tau}~(b)}
Fix $\varepsilon>0$. It follows from the
construction~(\ref{eq:def-XK}) of $\nu^K$ that, for $t<\tau_1$, under
$\mathbf{P}^K_{\frac{z_K}{K}\delta_x}$,
\begin{gather*}
  \nu^K_t=Z^K_t\delta_x, \\
  \mbox{where}\quad {\cal L}(Z^K)=\mathbf{P}^K((1-u_K\mu(x))b(x),
  d(x),\alpha(x,x),z_K/K). 
\end{gather*}
Therefore, by Theorem~\ref{thm:compar}~(c), for $K$ such that
$u_K<\varepsilon$ and for $t\leq\tau_1$,
\begin{gather}
  \label{eq:pf-lem-tau-(b)-1}
  Z^{K,1}_t\preceq\langle\nu^K_t,\mathbf{1}\rangle\preceq Z^{K,2}_t, \\
  \mbox{where}\quad
  \mathcal{L}(Z^{K,1})=\mathbf{P}^K((1-\varepsilon)b(x), 
  d(x),\alpha(x,x),z_K/K) \notag \\ \mbox{and}\quad
  \mathcal{L}(Z^{K,2})=\mathbf{P}^K(b(x),d(x),\alpha(x,x),z_K/K). \notag
\end{gather}

Now, let $\phi^1_y$, resp.~$\phi^2_y$, be the solution to
\begin{gather*}
  \dot{\phi}=((1-\varepsilon)b(x)-d(x)-\alpha(x,x)\phi)\phi, \\
  \mbox{resp.}\quad \dot{\phi}=(b(x)-d(x)-\alpha(x,x)\phi)\phi,
\end{gather*}
with initial state $y$, and observe that, for any $y>0$, when
$t\rightarrow+\infty$, $\phi^1_y(t)\rightarrow
e^1:=\bar{n}_x-\varepsilon b(x)/\alpha(x,x)$ and
$\phi^2_y(t)\rightarrow e^2:=\bar{n}_x$.

Define, for any $y>0$, $t^{i,y}_{\varepsilon}$ the first time such
that $\forall s\geq t^{i,y}_{\varepsilon}$,
$\phi^i_y(s)\in[e^i-\varepsilon,e^i+\varepsilon]$ ($i=1,2$). Because
of the continuity of the flows of these ODEs,
\begin{equation*}
  t^i_{\varepsilon}:=\sup_{y\in[z/2,2z]}t^{i,y}_{\varepsilon}<+\infty.
\end{equation*}

Let us apply Theorem~\ref{thm:exit-time}~(a) to $Z^{K,1}$ and
$Z^{K,2}$ on $[0,t_{\varepsilon}]$, where
$t_{\varepsilon}=t^1_{\varepsilon}\vee t^2_{\varepsilon}$: since
$z_K/K\rightarrow z$, for sufficiently small $\delta>0$, and for
$i=1,2$,
\begin{equation*}
  \lim_{K\rightarrow +\infty}\mathbf{P}
  \biggl(\sup_{0\leq t\leq t_{\varepsilon}}
  |Z^{K,i}_t-\phi^i_{z_K/K}(t)|>\delta\biggr)=0.
\end{equation*}
If we choose $\delta<\varepsilon$, we obtain, for $i=1,2$,
\begin{equation*}
  \lim_{K\rightarrow +\infty}\mathbf{P}
  (|Z^{K,i}_{t_{\varepsilon}}-e^i|<2\varepsilon)=1,
\end{equation*}
and so, for $i=1,2$,
\begin{equation}
  \label{eq:pf-lem-tau-(b)-2}
  \lim_{K\rightarrow +\infty}\mathbf{P}
  (|Z^{K,i}_{t_{\varepsilon}}-\bar{n}_x|<M\varepsilon)=1,
\end{equation}
where $M=2+b(x)/\alpha(x,x)$.

Now, assuming $\varepsilon$ sufficiently small for
$(M+1)\varepsilon<\bar{n}_x$, define the stopping times
\begin{equation*}
  T^{K,i}_{\varepsilon}=\inf\{t\geq
  t_{\varepsilon}:|Z^{K,i}_t-\bar{n}_x|>(M+1)\varepsilon\}
\end{equation*}
for $i=1,2$, and $T^K_{\varepsilon}=T^{K,1}_{\varepsilon}\wedge
T^{K,2}_{\varepsilon}$.

For any $z\in\mathbb{N}/K$, define also
\begin{equation*}
  \mathbf{P}_z^{K,1}:=\mathbf{P}^K((1-\varepsilon)b(x),d(x),\alpha(x,x),z).
\end{equation*}
Then, applying Theorem~\ref{thm:exit-time}~(c) to $\mathbf{P}_z^{K,1}$
with $C=[\bar{n}_x-M\varepsilon,\bar{n}_x+M\varepsilon]$, there exists
$V_1>0$ such that
\begin{gather}
  \lim_{K\rightarrow +\infty}\inf_{z\in
    C}\mathbf{P}^{K,1}_z(\hat{T}_{\varepsilon}>e^{KV_1})=1,
    \label{eq:th-3-explicite} \\
  \mbox{where}\quad \hat{T}_{\varepsilon}=\inf\{t\geq
    0:|w_t-\bar{n}_x|>(M+1)\varepsilon\}. \notag
\end{gather}
Therefore, applying the Markov property at time $t_{\varepsilon}$,
it follows from~(\ref{eq:pf-lem-tau-(b)-2}) that
\begin{equation*}
  \lim_{K\rightarrow +\infty}\mathbf{P}
  (T^{K,1}_{\varepsilon}>e^{KV_1}+t_{\varepsilon})=1.
\end{equation*}

Similarly, there exists $V_2>0$ such that
\begin{equation*}
  \lim_{K\rightarrow +\infty}\mathbf{P}
  (T^{K,2}_{\varepsilon}>e^{KV_2}+t_{\varepsilon})=1,
\end{equation*}
and thus
\begin{equation}
  \label{eq:pf-lem-tau-(b)-3}
  \lim_{K\rightarrow +\infty}\mathbf{P}
  (T^K_{\varepsilon}>e^{KV})=1,
\end{equation}
where $V:=V_1\wedge V_2$.

Now, because of~(\ref{eq:pf-lem-tau-(b)-1}),
\begin{equation}
  \label{eq:pf-lem-tau-(b)-4}
  \forall t\in[t_{\varepsilon},T^K_{\varepsilon}\wedge\tau_1],\quad
  |\langle\nu^K_s,\mathbf{1}\rangle-\bar{n}_x|<(M+1)\varepsilon.
\end{equation}
Therefore, since $\log K>t_{\varepsilon}$ for sufficiently large $K$,
in order to complete the proof of~(\ref{eq:lem-tau-(b)}), it suffices
to show that
\begin{equation}
  \label{eq:pf-lem-tau-(b)-5}
  \lim_{K\rightarrow +\infty}\mathbf{P}(\tau_1<T^K_{\varepsilon})=1.  
\end{equation}

If we denote by $A^K_t$ the number of mutations occuring between
$t_{\varepsilon}$ and $t+t_{\varepsilon}$, by
Theorem~\ref{thm:compar}~(b), for $t$ such that
$t_{\varepsilon}+t\leq T^K_{\varepsilon}\wedge\tau_1$,
\begin{equation*}
  B^K\preceq A^K,
\end{equation*}
where $B^K$ is a Poisson process with parameter
$Ku_K(\bar{n}_x-(M+1)\varepsilon)\mu(x)b(x)$.

Therefore, if we denote by $S^K$ the first time when $B^K_t=1$, on the
event $\{t_{\varepsilon}+S^K<T^K_{\varepsilon}\}$,
\begin{equation*}
  \tau_1\leq t_{\varepsilon}+S^K.  
\end{equation*}
Since $\exp(-KV)\ll Ku_K$,
$\lim_K\mathbf{P}(t_{\varepsilon}+S^K<e^{KV})=1$, and hence,
by~(\ref{eq:pf-lem-tau-(b)-3}),
\begin{equation*}
  \lim_{K\rightarrow +\infty}
  \mathbf{P}(t_{\varepsilon}+S^K<T^K_{\varepsilon})=1,
\end{equation*}
which implies~(\ref{eq:pf-lem-tau-(b)-5}).

In the case where $z_K/K\rightarrow\bar{n}_x$,
using~(\ref{eq:th-3-explicite}) as above, we obtain easily
\begin{gather*}
  \lim_{K\rightarrow +\infty}\mathbf{P}
  (S^K_{\varepsilon}>e^{KV})=1, \\
  \mbox{where}\quad S^K_{\varepsilon}=\inf\{t\geq
  0:|Z^{K,i}_t-\bar{n}_x|>(M+1)\varepsilon,\ i=1,2\}.
\end{gather*}
Then, the proof of~(\ref{eq:lem-tau-(b)-bar-n}) can be completed using
the same method as the one we used above.\hfill$\Box$

\paragraph{Proof of Lemma~\ref{lem:tau}~(c)}
Fix $t>0$ and $\varepsilon>0$. Take $K$ large enough for $\log
K<t/Ku_K$. The Markov property at time $\log K$ for $\nu^K$ yields
\begin{multline}
  \label{eq:pf-lem-tau-(c)-1}
  \mathbf{P}^K_{\frac{z_K}{K}\delta_x}\biggl(\tau_1>\frac{t}{Ku_K},\ 
  \sup_{t\in[\log K,\tau_1]}|\langle
  \nu^K_t,\mathbf{1}\rangle-\bar{n}_x|<\varepsilon\biggr) \\
  =\mathbf{E}^K_{\frac{z_K}{K}\delta_x}\biggl[\mathbf{1}_{\{\tau_1>\log
    K\}}\mathbf{P}^K_{\nu_{\log K}^K}
  \biggl(\tau_1>\frac{t}{Ku_K}-\log K,\ \\ \sup_{t\in[0,\tau_1]}|\langle
  \nu^K_t,\mathbf{1}\rangle-\bar{n}_x|<\varepsilon\biggr)\biggr].
\end{multline}

For any initial condition
$\nu^K_0=\langle\nu^K_0,\mathbf{1}\rangle\delta_x$ of $\nu^K$, by
Theorem~\ref{thm:compar}~(b), the number $A^K_t$ of mutations of
$\nu^K$ between 0 and $t$ satisfies, for any $t\leq\tau_1$ such that
$\sup_{s\in[0,t]}|\langle
\nu^K_s,\mathbf{1}\rangle-\bar{n}_x|<\varepsilon$,
\begin{equation*}
  B^K_t\preceq A^K_t\preceq C^K_t,
\end{equation*}
where $B^K_t$ and $C^K_t$ are Poisson processes with respective
parameters $Ku_K(\bar{n}_x-\varepsilon)\mu(x)b(x)$ and
$Ku_K(\bar{n}_x+\varepsilon)\mu(x)b(x)$.

Therefore, on the event $\{\sup_{s\in[0,\tau_1]}|\langle
\nu^K_s,\mathbf{1}\rangle-\bar{n}_x|<\varepsilon\}$,
$S^K\leq\tau_1\leq T^K$, where $T^K$ is the first time when $B^K_t=1$,
and $S^K$ the first time when $C^K_t=1$.

Now, by Lemma~\ref{lem:tau}~(b), under
$\mathbf{P}^K_{(z_K/K)\delta_x}$, $\nu^K_{\log K}\overset{\cal
  P}{\rightarrow}\bar{n}_x\delta_x$, so, by Skorohod's Theorem, we can
construct $\hat{N}^K$ with the same law as $\langle\nu^K_{\log
  K},\mathbf{1}\rangle$ on an auxiliary probability space
$\hat{\Omega}$ such that $\hat{N}^K(\hat{\omega})\rightarrow\bar{n}_x$
for any $\hat{\omega}\in\hat{\Omega}$. Fix
$\hat{\omega}\in\hat{\Omega}$. Then, by Lemma~\ref{lem:tau}~(b),
\begin{equation*}
  \lim_{K\rightarrow
    +\infty}\mathbf{P}^K_{\hat{N}(\hat{\omega})\delta_x}
  \biggl(\sup_{t\in[0,\tau_1]}|\langle
  \nu^K_t,\mathbf{1}\rangle-\bar{n}_x|<\varepsilon\biggr)=1,
\end{equation*}
and so,
\begin{multline*}
  \limsup_{K\rightarrow
    +\infty}\mathbf{P}^K_{\hat{N}(\hat{\omega})\delta_x}
  \biggl(\tau_1>\frac{t}{Ku_K}-\log K,\ \sup_{t\in[0,\tau_1]}|\langle
  \nu^K_t,\mathbf{1}\rangle-\bar{n}_x|<\varepsilon\biggr) \\
  \leq\limsup_{K\rightarrow
    +\infty}\mathbf{P}^K_{\hat{N}(\hat{\omega})\delta_x}
  \left(T^K>\frac{t}{Ku_K}-\log K\right)
  =\exp(-t(\bar{n}_x-\varepsilon)\mu(x)b(x)).
\end{multline*}
Therefore, under $\mathbf{P}^K_{(z_K/K)\delta_x}$,
\begin{multline*}
  \limsup_{K\rightarrow
    +\infty}\mathbf{P}^K_{\nu^K_{\log K}}
  \biggl(\tau_1>\frac{t}{Ku_K}-\log K,\ \sup_{t\in[0,\tau_1]}|\langle
  \nu^K_t,\mathbf{1}\rangle-\bar{n}_x|<\varepsilon\biggr) \\
  \leq\exp(-t(\bar{n}_x-\varepsilon)\mu(x)b(x))
\end{multline*}
in probability (where $\limsup X_n\leq a$ in probability means that,
for any $\eta>0$, $\mathbf{P}(X_n>a+\eta)\rightarrow 0$).

Similarly, under $\mathbf{P}^K_{(z_K/K)\delta_x}$,
\begin{multline*}
  \liminf_{K\rightarrow
    +\infty}\mathbf{P}^K_{\nu^K_{\log K}}
  \biggl(\tau_1>\frac{t}{Ku_K}-\log K,\ \sup_{t\in[0,\tau_1]}|\langle
  \nu^K_t,\mathbf{1}\rangle-\bar{n}_x|<\varepsilon\biggr) \\
  \geq\exp(-t(\bar{n}_x+\varepsilon)\mu(x)b(x))
\end{multline*}
in probability.

Now, by Lemma~\ref{lem:tau}~(a) and~(b),
\begin{gather*}
  \lim_{K\rightarrow +\infty}\mathbf{P}^K_{\frac{z_K}{K}\delta_x}
  (\tau_1>\log K)=1 \\
  \mbox{and}\quad
  \lim_{K\rightarrow +\infty}\mathbf{P}^K_{\frac{z_K}{K}\delta_x}
  \biggl(\sup_{t\in[\log K,\tau_1]}|\langle
  \nu^K_t,\mathbf{1}\rangle-\bar{n}_x|<\varepsilon\biggr)=1.
\end{gather*}
So, using property~(\ref{eq:cv-proba-bounded}), it follows
from~(\ref{eq:pf-lem-tau-(c)-1}) that
\begin{gather*}
  \limsup_{K\rightarrow
    +\infty}\mathbf{P}^K_{\frac{z_K}{K}\delta_x}
  \biggl(\tau_1>\frac{t}{Ku_K}\biggr)
  \leq\exp(-t(\bar{n}_x-\varepsilon)\mu(x)b(x)) \\
  \mbox{and}\quad \liminf_{K\rightarrow
    +\infty}\mathbf{P}^K_{\frac{z_K}{K}\delta_x}
  \biggl(\tau_1>\frac{t}{Ku_K}\biggr)
  \geq\exp(-t(\bar{n}_x+\varepsilon)\mu(x)b(x)).
\end{gather*}
Since this holds for any $\varepsilon>0$, we have completed the proof
of Lemma~\ref{lem:tau}~(c).\hfill$\Box$

\paragraph{Proof of Lemma~\ref{lem:theta}}
The proof of this lemma follows the three steps of the invasion of a
mutant described in Section~\ref{sec:outline} (cf.\ 
Fig.~\ref{fig:inv-fix}).

Fix $\eta>0$, $\varepsilon_0>0$ and $0<\varepsilon<\varepsilon_0$. By
Lemma~\ref{lem:tau}~(a), there exists a constant $\rho>0$ that we can
assume smaller than $\eta$, such that, for sufficiently large $K$,
\begin{equation}
  \label{eq:pf-lemma-theta-5}
  \mathbf{P}^K_{\frac{z_K}{K}\delta_x+\frac{1}{K}\delta_y}
  \left(\tau_1<\frac{\rho}{Ku_K}\right)<\varepsilon.
\end{equation}

Observe that, under
$\mathbf{P}^K_{\frac{z_K}{K}\delta_x+\frac{1}{K}\delta_y}$, for
$t\leq\tau_1$,
\begin{multline*}
  {\cal L}((\langle\nu^K,\mathbf{1}_{\{x\}}\rangle,
  \langle\nu^K,\mathbf{1}_{\{y\}}\rangle))
  =\mathbf{Q}^K((1-u_K\mu(x))b(x),(1-u_K\mu(y))b(y), \\ d(x),d(y),
  \alpha(x,x),\alpha(x,y),\alpha(y,x),\alpha(y,y),z_K/K,1/K).
\end{multline*}

Fix $K$ large enough for $u_K<\varepsilon$. Define
\begin{equation*}
  S^K_{\varepsilon}:=\inf\{s\geq 0:\langle
  \nu^K_s,\mathbf{1}_{\{y\}}\rangle\geq\varepsilon\}
\end{equation*}
By Theorem~\ref{thm:compar}~(c) and~(d), for $t<\tau_1\wedge
S^K_{\varepsilon}$,
\begin{gather}
  \label{eq:pf-lem-theta-2}
  Z^{K,1}_t\preceq\langle\nu^K_t,\mathbf{1}_{\{x\}}\rangle\preceq Z^{K,2}_t, \\
  \mbox{where}\quad {\cal L}(Z^{K,1})=\mathbf{P}^K
  ((1-\varepsilon)b(x),d(x)+\varepsilon\alpha(x,y),\alpha(x,x),z_K/K)
  \notag \\
  \mbox{and}\quad {\cal L}(Z^{K,2})=\mathbf{P}^K
  (b(x),d(x),\alpha(x,x),z_K/K). \notag
\end{gather}

Using the method that led us to~(\ref{eq:pf-lem-tau-(b)-3}), we can
deduce from Theorem~\ref{thm:exit-time}~(c) that there exists $V>0$
such that
\begin{gather}
  \label{eq:pf-lem-theta-3}
  \lim_{K\rightarrow +\infty}\mathbf{P}
  (R^K_{\varepsilon}>e^{KV})=1, \\
  \mbox{where}\quad R^K_{\varepsilon}=\inf\{t\geq
  0:|Z^{K,i}_t-\bar{n}_x|>M\varepsilon,\ i=1,2\}, \notag
\end{gather}
with $M=3+(b(x)+\alpha(x,y))/\alpha(x,x)$.

Now, observe that, by~(\ref{eq:pf-lem-theta-2}),
\begin{equation*}
  \forall t\leq\tau_1\wedge S^K_{\varepsilon}\wedge R^K_{\varepsilon},\quad
  \langle\nu^K_t,\mathbf{1}_{\{x\}}\rangle
  \in[\bar{n}_x-M\varepsilon,\bar{n}_x+M\varepsilon].
\end{equation*}
Therefore, by Theorem~\ref{thm:compar}~(c) and~(e), for
$t\leq\tau_1\wedge S^K_{\varepsilon}\wedge R^K_{\varepsilon}$
\begin{gather}
  \label{eq:pf-lem-theta-4}
  Z^{K,3}_t\preceq\langle\nu^K_t,\mathbf{1}_{\{y\}}\rangle\preceq
  Z^{K,4}_t, \quad\mbox{where} \\
  {\cal L}(Z^{K,3})=\mathbf{P}^K
  ((1-\varepsilon)b(y),d(y)+(\bar{n}_x+M\varepsilon)\alpha(y,x)
  +\varepsilon\alpha(y,y),0,1/K) \notag \\
  \mbox{and}\quad {\cal L}(Z^{K,4})=\mathbf{P}^K
  (b(y),d(y)+(\bar{n}_x-M\varepsilon)\alpha(y,x),0,1/K). \notag
\end{gather}

Define, for any $K\geq 1$, $n\in\mathbb{N}$ and $i\in\{3,4\}$, the
stopping time
\begin{equation*}
  T^{K,i}_{n/K}=\inf\{t\geq 0:Z^{K,i}_t=n/K\}.
\end{equation*}
Observe that, if $S^K_{\varepsilon}<\tau_1\wedge R^K_{\varepsilon}$,
\begin{equation}
  \label{eq:comp-Ti-S}
  T^{K,4}_{\lceil\varepsilon K\rceil/K}\leq
  S^K_{\varepsilon}\leq T^{K,3}_{\lceil\varepsilon K\rceil/K}
\end{equation}
and that, if $T^{K,4}_0<T^{K,4}_{\lceil\varepsilon
  K\rceil/K}\wedge\tau_1\wedge R^K_{\varepsilon}$,
\begin{equation*}
  \theta_0\leq T^{K,4}_0.
\end{equation*}

If $Z^{K,4}$ is sub-critical, apply
Theorem~\ref{thm:invasion-extinction}~(\ref{eq:non-inv-subcrit}), and
if $Z^{K,4}$ is super-critical, apply
Theorem~\ref{thm:invasion-extinction}~(\ref{eq:ext-supercrit}) (the
critical case can be excluded by slightly changing the value of
$\varepsilon$). Since $\log K\ll 1/Ku_K$, we obtain
\begin{multline}
  \label{eq:pf-lemma-theta-6}
  \lim_{K\rightarrow +\infty}
  \mathbf{P}\left(T^{K,4}_0\leq\frac{\rho}{Ku_K}\wedge
    T^{K,4}_{\lceil\varepsilon K\rceil/K}\right) \\ 
  =\frac{d(y)+(\bar{n}_x-M\varepsilon)\alpha(y,x)}{b(y)}\wedge 1
  \geq 1-\frac{[f(y,x)]_+}{b(y)}-\frac{\alpha(y,x)}{b(y)}M\varepsilon.
\end{multline}

Combining~(\ref{eq:pf-lemma-theta-5}), (\ref{eq:pf-lem-theta-3}),
(\ref{eq:pf-lem-theta-4}) and~(\ref{eq:pf-lemma-theta-6}), and using
the facts that $\rho<\eta$, $\varepsilon<\varepsilon_0$ and
$\exp(KV)>\rho/Ku_K$ for sufficiently large $K$, we obtain, taking $K$
larger if necessary,
\begin{align}
  \mathbf{P}
  \biggl(\theta_0<\tau_1\wedge & \frac{\eta}{Ku_K},\ V_0=x\mbox{\ and\ 
    }|\langle\nu^K_{\theta_0},\mathbf{1}\rangle-\bar{n}_x|<M\varepsilon_0
  \biggr) \notag \\
  & \geq\mathbf{P}\left(\theta_0<\tau_1\wedge S^K_{\varepsilon}\wedge
    R^K_{\varepsilon}\wedge\frac{\rho}{Ku_K}\mbox{\ and\ }V_0=x\right)
  \notag \\
  & \geq\mathbf{P}
  \left(T^{K,4}_0<\tau_1\wedge T^{K,4}_{\lceil\varepsilon
  K\rceil/K}\wedge
    R^K_{\varepsilon}\wedge\frac{\rho}{Ku_K}\right) \notag \\
  & \geq 1-\frac{[f(y,x)]_+}{b(y)}
  -\left(\frac{\alpha(y,x)}{b(y)}M+3\right)\varepsilon. \label{eq:pf-V0=x}  
\end{align}
This ends the proof of Lemma~\ref{lem:theta} in the case where
$f(y,x)\leq 0$.

Let us assume that $f(y,x)>0$, i.e.\ that
$b(y)-d(y)-\bar{n}_x\alpha(y,x)>0$. If we choose $\varepsilon>0$
sufficiently small, then $Z^{K,3}$ is super-critical. By
Theorem~\ref{thm:invasion-extinction}~(\ref{eq:inv-supercrit}),
\begin{align*}
  \lim_{K\rightarrow
  +\infty}\mathbf{P} & \left(T^{K,3}_{\lceil\varepsilon K\rceil/K}
  <\frac{\rho}{3Ku_K}\right) \\ & \qquad =\frac{(1-\varepsilon)b(y)-d(y)
    -(\bar{n}_x+M\varepsilon)\alpha(y,x)
    -\varepsilon\alpha(y,y)}{(1-\varepsilon)b(y)} \\
  & \qquad \geq\frac{f(y,x)}{(1-\varepsilon)b(y)}
  -\varepsilon\frac{b(y)+M\alpha(y,x)+\alpha(y,y)}{(1-\varepsilon)b(y)}.
\end{align*}
Therefore, by~(\ref{eq:pf-lem-theta-3})
and~(\ref{eq:pf-lemma-theta-5}), assuming (without loss of generality)
that $\varepsilon<1/2$, for sufficiently large $K$,
\begin{equation*}
  \mathbf{P}\left(T^{K,3}_{\lceil\varepsilon
      K\rceil/K}<\tau_1\wedge
      R^K_{\varepsilon}\wedge\frac{\rho}{3Ku_K}
  \right)\geq\frac{f(y,x)}{(1-\varepsilon)b(y)}-M'\varepsilon,
\end{equation*}
where $M':=2(b(y)+M\alpha(y,x)+\alpha(y,y))/b(y)+3$. Then, it follows
from~(\ref{eq:comp-Ti-S}) that
\begin{equation}
  \label{eq:pf-lemma-theta-10}
  \mathbf{P}\left(S^K_{\varepsilon}<\tau_1\wedge
    R^K_{\varepsilon}\wedge\frac{\rho}{3Ku_K} 
  \right)\geq\frac{f(y,x)}{(1-\varepsilon)b(y)}-M'\varepsilon.
\end{equation}

Observe that, on the event $\{S^K_{\varepsilon}<\tau_1\wedge
R^K_{\varepsilon}\wedge(\rho/3Ku_K)\}$,
\begin{equation}
  \label{eq:nuK-at-time-S}
  \langle\nu^K_{S^K_{\varepsilon}},\mathbf{1}_{\{y\}}\rangle
  =\lceil\varepsilon K\rceil/K \quad\mbox{and}\quad
  |\langle\nu^K_{S^K_{\varepsilon}},\mathbf{1}_{\{x\}}\rangle
  -\bar{n}_x|<M\varepsilon.
\end{equation}

Now, since we have assumed $f(y,x)>0$, $x$ and $y$
satisfy~(\ref{eq:hyp-B2}) and, by Proposition~\ref{prop:LV}, any
solution to~(\ref{eq:LV-IPS}) with initial state in the compact set
$[\bar{n}_x-M\varepsilon,\bar{n}_x+M\varepsilon]
\times[\varepsilon/2,2\varepsilon]$ converges to $(0,\bar{n}_y)$ when
$t\rightarrow +\infty$. As in the proof of Lemma~\ref{lem:tau} (b),
because of the continuity of the flow of system~(\ref{eq:LV-IPS}), we
can find $t_{\varepsilon}<+\infty$ large enough such that any of these
solutions do not leave the set
$[0,\varepsilon^2/2]\times[\bar{n}_y-\varepsilon/2,\bar{n}_y+\varepsilon/2]$
after time $t_{\varepsilon}$.

Apply Theorem~\ref{thm:exit-time}~(b) on $[0,t_{\varepsilon}]$, with
$C=[\bar{n}_x-M\varepsilon,\bar{n}_x+M\varepsilon]
\times[\varepsilon/2,2\varepsilon]$ and with a constant
$\delta<\varepsilon^2/2\wedge r$, where $r$ is defined
in~(\ref{eq:def-r}) (with $T=t_{\varepsilon}$). Then, with the
notations of Theorem~\ref{thm:exit-time}~(b), because
of~(\ref{eq:pf-lemma-theta-10}) and~(\ref{eq:nuK-at-time-S}), the
Markov property at time $S^K_{\varepsilon}$ yields
\begin{multline}
  \label{eq:pf-lemma-theta-11}
  \liminf_{K\rightarrow +\infty}\mathbf{P}
  \biggl(S^K_{\varepsilon}<\tau_1\wedge
    R^K_{\varepsilon}\wedge\frac{\rho}{3Ku_K},\  \\
    \sup_{S^K_{\varepsilon}\leq s\leq
    S^K_{\varepsilon}+t_{\varepsilon}}
  \bigl\|\bigl(\langle\nu^K_s,\mathbf{1}_{\{x\}}\rangle,
  \langle\nu^K_s,\mathbf{1}_{\{y\}}\rangle\bigr)
  -\phi_{\langle\nu^K_{S^K_{\varepsilon}},\mathbf{1}_{\{x\}}\rangle,
  \langle\nu^K_{S^K_{\varepsilon}},\mathbf{1}_{\{y\}}\rangle}(s)\bigr\|
      \leq\delta\biggr) \\
      \geq\frac{f(y,x)}{(1-\varepsilon)b(y)}-M'\varepsilon.
\end{multline}
Now, observe that, since $\delta<r$, on the event
\begin{multline*}
  \biggl\{S^K_{\varepsilon}<\tau_1\wedge
    R^K_{\varepsilon},\ \\
    \sup_{S^K_{\varepsilon}\leq s\leq
    S^K_{\varepsilon}+t_{\varepsilon}}
  \bigl\|\bigl(\langle\nu^K_s,\mathbf{1}_{\{x\}}\rangle,
  \langle\nu^K_s,\mathbf{1}_{\{y\}}\rangle\bigr)
  -\phi_{\langle\nu^K_{S^K_{\varepsilon}},\mathbf{1}_{\{x\}}\rangle,
  \langle\nu^K_{S^K_{\varepsilon}},\mathbf{1}_{\{y\}}\rangle}(s)\bigr\|
      \leq\delta\biggr\},
\end{multline*}
for any $t\in[S^K_{\varepsilon},S^K_{\varepsilon}+t_{\varepsilon}]$,
$\langle\nu^K_t,\mathbf{1}_{\{x\}}\rangle\geq r-\delta>0$ and
$\langle\nu^K_t,\mathbf{1}_{\{y\}}\rangle\geq r-\delta>0$, and thus
\begin{equation*}
  \theta_0>S^K_{\varepsilon}+t_{\varepsilon}.
\end{equation*}

Therefore, since $\delta<\varepsilon^2/2<\varepsilon/2$,
by~(\ref{eq:pf-lemma-theta-5}) and~(\ref{eq:pf-lemma-theta-11}), for
sufficiently large $K$,
\begin{multline}
  \label{eq:pf-lemma-theta-12}
  \mathbf{P}\left(S^K_{\varepsilon}<R^K_{\varepsilon}
    \wedge\frac{\rho}{3Ku_K},\ \tau_1>\frac{\rho}{3Ku_K}+t_{\varepsilon},\ 
    \theta_0>S^K_{\varepsilon}+t_{\varepsilon},\ \right. \\ \left. \langle
    \nu^K_{S^K_{\varepsilon}+t_{\varepsilon}},
    \mathbf{1}_{\{x\}}\rangle<\varepsilon^2\mbox{\ and\ }
    \langle\nu^K_{S^K_{\varepsilon}+t_{\varepsilon}},
    \mathbf{1}_{\{y\}}\rangle
    \in[\bar{n}_y-\varepsilon,\bar{n}_y+\varepsilon]\right) \\ \geq
  \frac{f(y,x)}{(1-\varepsilon)b(y)}-(M'+2)\varepsilon.
\end{multline}

Now, we will compare $\langle\nu^K,\mathbf{1}_{\{x\}}\rangle$ with a
branching process after time $S^K_{\varepsilon}+t_{\varepsilon}$ in
order to prove that trait $x$ gets extinct with a very high
probability. We will use a method very similar to the one we used in
the beginning of this proof. First, on the event inside the
probability in~(\ref{eq:pf-lemma-theta-12}),
$\langle\nu^K_{S^K_{\varepsilon}+t_{\varepsilon}},
\mathbf{1}_{\{x\}}\rangle<\varepsilon^2$. In order to prove that the
population with trait $x$ stays small after
$S^K_{\varepsilon}+t_{\varepsilon}$, let us define the stopping time
\begin{equation*}
  \hat{S}^K_{\varepsilon}=\inf\{t\geq
  S^K_{\varepsilon}+t_{\varepsilon}:\langle
  \nu^K_t,\mathbf{1}_{\{x\}}\rangle>\varepsilon\}
\end{equation*}
(remind that $\varepsilon^2<\varepsilon$ since $\varepsilon<1/2$).
Using Theorem~\ref{thm:compar}~(c) and~(d) again, we see that, on the
event
\begin{equation*}
  F^{K,\varepsilon}:=\bigl\{\langle
  \nu^K_{S^K_{\varepsilon}+t_{\varepsilon}},
  \mathbf{1}_{\{x\}}\rangle<\varepsilon^2,\ 
  \langle\nu^K_{S^K_{\varepsilon}+t_{\varepsilon}},
  \mathbf{1}_{\{y\}}\rangle
  \in[\bar{n}_y-\varepsilon,\bar{n}_y+\varepsilon]\bigr\},
\end{equation*}
for any $t\geq 0$ such that
$S^K_{\varepsilon}+t_{\varepsilon}+t\leq\hat{S}^K_{\varepsilon}\wedge\tau_1$,
\begin{gather*}
  Z^{K,5}_t\preceq\langle\nu^K_{S^K_{\varepsilon}+t_{\varepsilon}+t},
  \mathbf{1}_{\{y\}}\rangle\preceq Z^{K,6}_t, \\
  \mbox{where}\quad {\cal L}(Z^{K,5})=\mathbf{P}^K((1-\varepsilon)b(y),
  d(y)+\varepsilon\alpha(y,x),\alpha(y,y),
  \lfloor(\bar{n}_y-\varepsilon)K\rfloor/K) \\
  \mbox{and}\quad {\cal L}(Z^{K,6})=\mathbf{P}^K(b(y),d(y),
  \alpha(y,y),\lceil(\bar{n}_y+\varepsilon)K\rceil/K).
\end{gather*}

We can apply Theorem~\ref{thm:exit-time}~(c) to $Z^{K,5}$ and
$Z^{K,6}$ as above to obtain a constant $V'>0$ such that
\begin{gather}
  \label{eq:pf-lem-theta-13}
  \lim_{K\rightarrow+\infty}\mathbf{P}(\hat{R}^K_{\varepsilon}>e^{KV'})=1, \\
  \mbox{where}\quad \hat{R}^K_{\varepsilon}=
  \inf\{t\geq 0:|Z^{K,i}_t-\bar{n}_y|>M''\varepsilon, i=5,6\}, \notag
\end{gather}
with $M''=3+(b(y)+\alpha(y,x))/\alpha(y,y)$.

Observe that, on the event $F^{K,\varepsilon}$, for any
$t\leq\hat{R}^K_{\varepsilon}$ such that
$S^K_{\varepsilon}+t_{\varepsilon}+t
\leq\hat{S}^K_{\varepsilon}\wedge\tau_1$,
\begin{equation*}
  |\langle\nu^K_{S^K_{\varepsilon}+t_{\varepsilon}+t},
  \mathbf{1}_{\{y\}}\rangle-\bar{n}_y|\leq M''\varepsilon,
\end{equation*}
and so, by Theorem~\ref{thm:compar}~(c) and~(e), on
$F^{K,\varepsilon}$ and for $t$ as above,
\begin{gather*}
  \langle\nu^K_{S^K_{\varepsilon}+t_{\varepsilon}+t},
  \mathbf{1}_{\{x\}}\rangle\preceq Z^{K,7}_t \\
  \mbox{where}\quad {\cal L}(Z^{K,7})=\mathbf{P}^K(b(x),
  d(x)+(\bar{n}_y-M''\varepsilon)\alpha(x,y),0,\lceil\varepsilon^2K\rceil/K).
\end{gather*}

Now, since $x$ and $y$ satisfy~(\ref{eq:hyp-B2}), $Z^{K,7}$ is
sub-critical for sufficiently small $\varepsilon$. Fix such an
$\varepsilon>0$ and define for any $n\geq 0$
\begin{equation*}
  \hat{T}^K_{n/K}=\inf\{t\geq 0:Z^{K,7}_t=n/K\}.
\end{equation*}
If $\hat{T}^K_{\lceil\varepsilon
  K\rceil/K}\leq\hat{R}^K_{\varepsilon}$ and
$S^K_{\varepsilon}+t_{\varepsilon}+\hat{T}^K_{\lceil\varepsilon
  K\rceil/K}\leq\tau_1$, then
\begin{equation*}
  \hat{S}^K_{\varepsilon}\geq S^K_{\varepsilon}+t_{\varepsilon}
  +\hat{T}^K_{\lceil\varepsilon K\rceil/K}
\end{equation*}
and if $\hat{T}^K_0\leq\hat{R}^K_{\varepsilon}$ and
$S^K_{\varepsilon}+t_{\varepsilon}+\hat{T}^K_0\leq
\hat{S}^K_{\varepsilon}\wedge\tau_1$, then
\begin{equation*}
  \theta_0\leq\hat{T}^K_0.
\end{equation*}

Moreover, by
Theorem~\ref{thm:invasion-extinction}~(\ref{eq:ext-subcrit})
and~(\ref{eq:ext-before-expl}), for sufficiently large $K$,
\begin{gather*}
  \mathbf{P}\left(\hat{T}^K_0
    \leq\frac{\rho}{3Ku_K}\right)\geq 1-\varepsilon \\
  \mbox{and}\quad \mathbf{P}(\hat{T}^K_{\lceil
    K\varepsilon\rceil/K}\leq\hat{T}^K_0)\leq 2\varepsilon.
\end{gather*}

Combining the last two inequalities with~(\ref{eq:pf-lemma-theta-5}),
(\ref{eq:pf-lemma-theta-12}) and (\ref{eq:pf-lem-theta-13}), and
reminding that $\rho<\eta$ and $\varepsilon<\varepsilon_0$, we finally
obtain, for sufficiently large $K$,
\begin{align*}
  & \mathbf{P}\left(\theta_0<\tau_1\wedge\frac{\eta}{Ku_K},\ 
    V_0=y\mbox{\ and\ 
    }|\langle\nu^K_{\theta_0},\mathbf{1}\rangle-\bar{n}_y|
    <M''\varepsilon_0\right) \\
  & \!\geq\mathbf{P}
  \left(S^K_{\varepsilon}<R^K_{\varepsilon}\wedge\frac{\rho}{3Ku_K},\ 
    \!\theta_0>S^K_{\varepsilon}+t_{\varepsilon},\ \!
    \tau_1>\frac{2\rho}{3Ku_K}+t_{\varepsilon},\ \!
    \langle\nu^K_{S^K_{\varepsilon}+t_{\varepsilon}},
    \mathbf{1}_{\{x\}}\rangle<\varepsilon^2,\ \! \right. \\ &
  \quad\left.
    \langle\nu^K_{S^K_{\varepsilon}+t_{\varepsilon}},\mathbf{1}_{\{y\}}
    \rangle\in[\bar{n}_y-\varepsilon,\bar{n}_y+\varepsilon],\ 
    \hat{T}^K_0<\frac{\rho}{3Ku_K}\wedge\hat{T}^K_{\lceil
      K\varepsilon\rceil/K}\mbox{\ and\ 
    }\hat{R}^K_{\varepsilon}>\frac{\rho}{Ku_K}\right) \\
  & \!\geq\frac{f(y,x)}{(1-\varepsilon)b(y)}-(M'+7)\varepsilon.
\end{align*}
Adding this inequality with~(\ref{eq:pf-V0=x}), we obtain
\begin{equation*}
  \mathbf{P}\left(\theta_0<\tau_1\wedge\frac{\eta}{Ku_K}\right)\geq
  1-\frac{\varepsilon}{1-\varepsilon}
  \frac{f(y,x)}{b(y)}-\left(M\frac{\alpha(y,x)}{b(y)}+M'+10\right)
  \varepsilon\geq 1-M'''\varepsilon,
\end{equation*}
where $M'''=2f(y,x)/b(y)+M\alpha(y,x)/b(y)+M'+10$, which
implies~(\ref{eq:lem-theta-3}), and
\begin{equation*}
  \mathbf{P}\bigl(|\langle\nu^K_{\theta_0},\mathbf{1}\rangle
  -\bar{n}_{V_0}|<(M\vee
  M'')\varepsilon_0\bigr)\geq 1-M'''\varepsilon,
\end{equation*}
which implies~(\ref{eq:lem-theta-4}).

Therefore,
\begin{equation*}
  \mathbf{P}(V_0=x)\geq 1-\frac{f(y,x)}{b(y)}-2M'''\varepsilon
  \quad\mbox{and}\quad
  \mathbf{P}(V_0=y)\geq\frac{f(y,x)}{(1-\varepsilon)b(y)}-2M'''\varepsilon.
\end{equation*}
Since $\mathbf{P}(V_0=x)\leq1-\mathbf{P}(V_0=y)$, we finally
obtain~(\ref{eq:lem-theta-1}) and~(\ref{eq:lem-theta-2}).\hfill$\Box$


\bigskip

\textbf{Acknowledgments:} I would like to thank S. M\'el\'eard and R.
Ferri\`ere for their continual guidance during my work. I also thank
the referees for useful comments and suggestions.

\end{document}